\documentclass[12pt,reqno]{amsart}
\usepackage{a4wide}
\allowdisplaybreaks \numberwithin{equation}{section}
\usepackage{color}

\numberwithin{equation}{section}

\newtheorem{theorem}{Theorem}[section]
\newtheorem{proposition}[theorem]{Proposition}

\newtheorem{lemma}[theorem]{Lemma}

\theoremstyle{definition}

\theoremstyle{remark}
\newtheorem{remark}[theorem]{Remark}

\newcommand{\ep}{\varepsilon}
\newcommand{\Om}{\Omega}

\newcommand{\ds}{\displaystyle}

\begin{document}

\title[uniqueness of  planar vortex patch
 ]
{uniqueness of  planar vortex patch  in incompressible steady flow }

 \author{Daomin Cao, Yuxia Guo, Shuangjie Peng   and Shusen Yan
}

\address{Institute of Applied Mathematics, Chinese Academy of Science, Beijing 100190, and University of Chinese Academy of Sciences, Beijing 100049,  P.R. China}
\email{dmcao@amt.ac.cn}

\address{  Department of Mathematical Science, Tsinghua University, Beijing, P.R.China}
\email{yguo@math.tsinghua.edu.cn}

\address{School of  Mathematics and Statistics \& Hubei Key Laboratory of Mathematical Sciences, Central China Normal University, Wuhan,
 P.R. China}

\email{
sjpeng@mail.ccnu.edu.cn}

\address{ Department of Mathematics, The  University of New England Armidale, NSW 2351, Australia}

\email{
syan@turing.une.edu.au
}


\begin{abstract}

We investigate  a
steady planar flow of an ideal fluid in a bounded simple connected  domain  and focus
on the   vortex patch problem with prescribed vorticity
strength. There are two methods to deal with the existence of solutions for this problem: the vorticity method and the stream function method.
A long standing open problem is whether these two entirely different methods result in the same solution.
In this paper, we will give a positive answer to this problem by studying the local uniqueness of the solutions. Another result obtained in this paper is that
if the domain is convex, then the vortex patch problem  has a unique solution.

\end{abstract}

\maketitle

\section{Introduction}

The incompressible steady flow without external  force is governed by the following  mass equation
\begin{equation}\label{1-4-11}
\nabla\cdot\mathbf{v}=0,
\end{equation}
and the following    Euler motion equations
\begin{equation}\label{2-4-11}
(\mathbf{v}\cdot\nabla)\mathbf{v}=-\nabla P,
\end{equation}
where $\mathbf{v}$ is  the velocity  and $P$ is the pressure
 in  the flow. Here, we assume the density is one.  If we consider a flow in a domain $\Omega$, we usually impose the
 following impermeable boundary condition:
 \begin{equation}\label{1-7-9}
\mathbf{v}\cdot \nu=0,
\end{equation}
 where $\nu$ is the outward unit normal of $\partial \Omega$.

  Introducing the vorticity vector $\vec{\omega} = curl \mathbf v$, we can rewrite the Euler equation as
\begin{equation}\label{n2-4-11}
\vec{\omega}\times \mathbf{v}=-\nabla \bigl(P+ \frac12 |\mathbf v|^2\bigr).
\end{equation}

 In this paper, we will consider the planar flow in the domain $\Omega$ in $\mathbb R^2$.  So $  \mathbf{v}= (v_1, v_2, 0)$, and     $\vec \omega=\bigl(0, 0,\omega\bigr)$, where
$\omega=  \frac{\partial
v_2}{\partial x_1}-\frac{\partial v_1}{\partial x_2}$.  On the other hand,
it follows from \eqref{1-4-11} and \eqref{1-7-9} that  for an incompressible steady planar flow, in any  connected
 domain $\Omega$ (not necessarily simple connected), there is a function $\psi$, which
is called the  stream function of the flow,  such that
\begin{equation}\label{3-4-11}
 \mathbf{v}= \bigl(\frac{\partial \psi}{\partial x_2},\; -\frac{\partial \psi}{\partial x_1}\bigr), \quad \text{in}\; \Omega.
\end{equation}
Then the vorticity can be written as
\begin{equation}\label{4-4-11}
 \omega=\partial_1
v_2-\partial_2 v_1=-\Delta \psi.
\end{equation}
and
\begin{equation}\label{2-7-9}
\vec{\omega}\times \mathbf{v}= \omega \nabla \psi.
\end{equation}
Moreover, \eqref{1-7-9} implies on each connected component of $\partial \Omega$, $\psi$ is a constant.  So, if $\Omega$ is simple connected,
\begin{equation}\label{3-7-9}
\psi(x)=0,\quad x\in\partial\Omega,
\end{equation}
after suitably adding a constant to $\psi$. In this paper, we always assume that $\Omega$ is simple connected.

The question on the existence of solutions representing steady
vortex rings occupies a central place in the theory of vortex motion
initiated by Helmholtz in 1858. See for example \cite{AS,AF,AT,BF,B1,B2,CLW,CPY1,FB,FB1,FT,L,Ni,N,N1,SV,T,T2} and the references therein.
In this paper, we will consider a
steady planar flow of an ideal fluid in a bounded region and focus
on the    flow, whose
vorticity $\omega$ is a constant $ \lambda$ in a region
$\Omega_\lambda$ which has $k$ connected components $
\Omega_{\lambda, j}$ and $ \Omega_{\lambda, j}$ shrinks to $k$
points  $x_{0, j}\in\bar \Omega$,
 $j=1, \cdots, k$, as $\lambda\to+\infty$, while  $\omega=0$ elsewhere. Moreover, it holds
\begin{equation}\label{3-7-9}
\int_{\Omega_{\lambda, j}} \omega =\kappa_j,
\end{equation}
where $\kappa_j>0$ is a given constant.
Such problem is called the   vortex patch problem with prescribed vorticity
strength at each vortex point.
Here, we do not assume that
$x_{0, i}\ne x_{0, j}$ for $i\ne j$, nor  $x_{0, j}\in \Omega$.

 Write
\[
\Omega_\lambda= \cup_{j=1}^k \Omega_{\lambda, j}.
\]

 From the above discussion, we find that $\psi$ satisfies the following elliptic problem:
\begin{equation}\label{4-7-9}
\begin{cases}
- \Delta \psi =\lambda \sum\limits_{j=1}^k  1_{\Omega_{\lambda, j}} ,&\text{in}\;\Om,\\
\psi=0, &\text{on}\; \partial\Om,
\end{cases}
\end{equation}
where $1_S=1$ in $S$ and $1_S=0$ elsewhere for any non-empty set  $S$.

 Let us point out that the Euler equation will give a relation
 between the set $\Omega_{\lambda, j}$ and $\psi$.
 Indeed, it follows from \eqref{n2-4-11}
and \eqref{2-7-9}  that the following relation holds,
\begin{equation}\label{5-7-9}
\lambda \sum\limits_{j=1}^k  1_{\Omega_{\lambda, j}}\nabla \psi=-\nabla \bigl(P+ \frac12 |\mathbf v|^2\bigr),
\end{equation}
which implies that $\Omega_{\lambda, j} =B_\delta(x_{0, j})
 \cap\{\psi>\tilde \kappa_{\lambda,j}\}$, where $x_{0,j}$ is the point that $\Omega_{\lambda,j}$ is assumed to shrink to as $\lambda\to+\infty$. See the discussion in Lemma~\ref{l1-11-9}.  Thus, \eqref{4-7-9} becomes
\begin{equation}\label{7-4-11}
\begin{cases}
- \Delta \psi =\lambda \sum\limits_{j=1}^k 1_{B_\delta(x_{0, j})}
1_{ \{\psi>\tilde \kappa_{\lambda,j}\}},&\text{in}\;\Om,\\
\psi=0, &\text{on}\; \partial\Om,
\end{cases}
\end{equation}
for some large  $\tilde \kappa_{\lambda,j}\ge \kappa>0$, subject to the following prescribed vortex strength condition
 \begin{equation}\label{10-7-9}
\lambda |\Omega_{\lambda, j}| =\kappa_j>0.
\end{equation}

Let us remark that once we find the stream function $\psi$, the velocity of the flow is given
by \eqref{3-4-11} and the pressure is given by
\begin{equation}\label{n5-4-11}
 P=\lambda \sum_{j=1}^k 1_{B_\delta(x_{0, j})} (\psi -\tilde \kappa_{\lambda, j})_+-\frac12 |\nabla\psi|^2,
\end{equation}
where
$\psi_+ =\psi$ if $\psi\ge 0$  and $\psi_+=0$ if $\psi<0$.

Let
$G$ be the Green function for $-\Delta$ in $\Omega$ with zero
boundary condition, written as
\[
G(x,x')=\frac{1}{2\pi}\ln \frac1{|x-x'|}-H(x,x'),   \,\,\,x,x'\in
\Omega.
\]
Recall that the Robin function is defined as
\[
\varphi(x)=H(x,x).
\]

 For any given integer $k>0$, we define the following Kirchhoff-Routh function (see \cite{L}):
\begin{equation}\label{KR}
\mathcal{W}(x_1,\cdots,x_k)=-\sum_{i\neq j}^k\kappa_i\kappa_j
G(x_i,x_j)+\sum^{k}_{i=1}\kappa_i^2 \varphi(x_i).
\end{equation}
Note that if $k=1$, then $\mathcal{W}=\kappa^2 \varphi$.

In \cite{CPY1}, we prove the following existence result:

\medskip

\textbf{Theorem~A.}  {\it
Suppose that ${\bf  x}_0 \in\Omega^k$ is an isolated critical point of
$\mathcal{W}({\bf  x})$,  satisfying $ \deg(\nabla \mathcal{W}, {\bf  x}_0)\ne 0$. Then, there
is an $\lambda_0>0$, such that for all $\lambda\in (\lambda_0, +\infty)$, \eqref{7-4-11}--\eqref{10-7-9}  has
a solution $\psi_\lambda$ such that the vorticity set $\{ x:   \; \omega_\lambda(x)=\lambda \} $ shrinks
to ${\bf  x}_0$ as $\lambda\to +\infty$.
}

\medskip

Equation~\eqref{7-4-11} has jumping nonlinearities. Its solutions are not in $C^2$.  This kind of discontinuity are much more
difficult to deal with than those in \cite{CF,CPY}, where the derivative of the nonlinearity is discontinuous.
  The proof of Theorem~A involves
the domain variation type estimates.

In \cite{T}, Turkington considered  the vortex patch problem with prescribed vorticity
strength \eqref{7-4-11}--\eqref{10-7-9} for the case $k=1$.
He obtained an existence result by studying the asymptotic behavior of the   absolute maximizer of   the kinetic
energy defined by
\[
\frac12\int_\Omega |\mathbf v|^2dx=\frac{1}{2}\int_\Omega \int_\Omega
\omega(x)G(x,y)\omega(y)dxdy
\]
in the following class
\[
K_\lambda(\Omega)=\left\{\omega\in
L^\infty(\Omega)\,:\,\int_\Omega\omega(x)dx=\kappa,
0\leq\omega(x)\leq\lambda \,\,\,a.e.\,\, x\in \Omega\right\}.
\]
 It was   proved in \cite{T} that the maximizer $\omega_\lambda$ satisfies
 \begin{equation}\label{add5}
 \omega_\lambda =\lambda 1_{\Omega_\lambda}, \Omega_\lambda=\{x\in \Omega \,:\,\psi_\lambda(x)>0\},
 \end{equation}
  where
$\psi_\lambda$ is the corresponding stream function satisfying
\begin{equation}\label{psi}
\begin{cases}
 -\Delta \psi=\lambda 1_{\{\psi(x)>0\}},\,\,\text{in}\,\Omega,\\
\psi=\mu_\lambda,\,\, \text{on}\,\partial\Omega,
\end{cases}
\end{equation}
for some constant  $\mu_\lambda$,    which satisfies
$\mu_\lambda =-\log\lambda+O(1)$ for $\lambda$ large.
Moreover, $\Omega_\lambda$ shrinks to a point $x_0$, which is a global minimum point of the Robin function $\varphi(x)$.

 Based on  the above mentioned results, Turkington pointed out that the geometry of $\Omega$ may lead to
 the non-uniqueness of solutions for \eqref{psi}. Theorem~A
 confirms this observation by establishing a relation between the existence of solutions for  \eqref{psi} and the  non-degenerate critical points of
 the Robin function $\varphi$. To prove Theorem~A, we work on the stream function $\psi$ instead of the vorticity function $\omega$.
 The stream function method has the advantage to obtain solutions with the vorticity set shrinking to the saddle point of the
 function $\mathcal{W}({\bf  x})$, while the vorticity method has strong physical motivation.
 It was asked in \cite{FT} whether these two entirely different methods give the same solutions. This is a local uniqueness
 problem.
 However, as far as we know, no much is known on the uniqueness of solution for  the vortex patch problem with prescribed vorticity
strength \eqref{7-4-11}--\eqref{10-7-9}.  The aim of this paper is to study the local uniqueness of solution for this problem and thus
prove that the vorticity method and the stream function method just give the same solution.

One of the main results of this paper is the following.

\begin{theorem}\label{th11}

Let $\kappa_j$, $j=1,\cdots,k$,  be $k$ given positive numbers.
Suppose that $\psi_\lambda$ is a solution of \eqref{7-4-11} and
\eqref{10-7-9},  such that each component of  vorticity set $\Omega_{\lambda,j}\,(j=1,\cdots,k)$ shrinks to $x_{0, j}\in \bar\Omega$,  as
$\lambda\to +\infty$. Then, $x_{0, j}\in \Omega$,  $j=1,\cdots, k$,  $x_{0, j}\ne x_{0, i}$ for $j\ne i$, and
$ {\bf
x}_0=(x_{0,1}, \cdots, x_{0, k})$ is a critical point of $\mathcal{W}$.

\end{theorem}

Theorem~A and Theorem~\ref{th11} show that the existence of solutions for \eqref{7-4-11} and
\eqref{10-7-9} are nearly determined by the critical points of the function $\mathcal{W}$. Results on the existence
and non-degeneracy of critical points for $\mathcal{W}$ can be found in \cite{BP,BPW}.
In  \cite{GT}, it was proved that there  does not exist  any critical point of $\mathcal{W}(x_1,\cdots,x_k)$ in $\Omega^k$ with $k\ge 2$ and $\kappa_j>0$ for all $j=1,\cdots,k$ if $\Omega$ is convex.
Hence, a direct consequence of Theorem~\ref{th11} is that if $\Omega$ is a convex domain,  \eqref{7-4-11}--\eqref{10-7-9} has no solution for $k\ge 2$. To obtain
a uniqueness result in convex domains, we  need to consider the case $k=1$.
  Another main result
in this paper is the following uniqueness result.

\begin{theorem}\label{th1-21-10}

Suppose that $x_0$ is an isolated  critical point of $\varphi(x)= H(x, x)$, which is  non-degenerate. Then  for large $\lambda>0$,
\eqref{7-4-11} with $k=1$, together with
\eqref{10-7-9},   has a unique solution.
\end{theorem}

The local uniqueness result in Theorem~\ref{th1-21-10} shows that any non-degenerate critical point of the Robin function $\varphi$ can only generate one solution
for \eqref{7-4-11} with $k=1$. This result clearly implies that the vorticity method and the stream function method actually result in the same solution.
On the other hand,
if $k=1$, then it follows from \cite{CF1} that in a convex domain,  $\mathcal{W}=\kappa^2 \varphi$ has a unique critical point, which is also non-degenerate. This result and
Theorem~\ref{th1-21-10} give the following uniqueness result in convex domains.

\begin{theorem}\label{th10}

Suppose that $\Omega$ is convex. Then the vortex patch problem with prescribed vorticity
strength \eqref{7-4-11}--\eqref{10-7-9} has a unique solution  and $k=1$ if $\lambda>0$ is large.
\end{theorem}

Our uniqueness result shows that if the domain is convex, then the flow can only has one vortex and the vortex point must be near the unique global
minimum point of the Robin function  $\varphi(x)$. Moreover, the vorticity of this solution must be the maximizer of the kinetic energy which was studied  by Turkington in \cite{T}.

The paper is organized as follows. In sections~2 and 3, we will study the asymptotic behavior of the solutions. Theorems~\ref{th11}
and \ref{th10} are proved in section~4 and section~5 respectively. The discussion of the free boundary $\partial\Omega_{\lambda, j}$ is given
in the appendix.

To analyze  the asymptotic behavior of the solutions for \eqref{7-4-11}, it is important to determine the scalar in the
blow-up procedure and find the corresponding limit problem. Obviously,  it is more reasonable to scale the equation in  \eqref{7-4-11} by using the diameter $D_{\lambda, j}$ of the unknown set
  $\Omega_{\lambda, j}$. Note that in  \eqref{7-4-11}, the parameter $\tilde \kappa_{\lambda, j}$ is also unknown. So,
  the crucial step is  to estimate both $D_{\lambda, j}$ and $\tilde \kappa_{\lambda, j}$ in terms of $\lambda$ and $\kappa_j$.
  These are achieved by using the Pohozaev identity and Harnack inequality.
  Let us point out that the estimates for \eqref{7-4-11} are domain variation type estimates in view of the terms
$1_{ \{\psi>\tilde \kappa_{\lambda,j}\}}$ appearing in \eqref{7-4-11}. Once we obtain the asymptotic of the solutions $\psi_\lambda$,
we can use the Pohozaev identity to prove Theorem~\ref{th11}.

The discussion of the local uniqueness of concentration solutions is dated back to the early 1990s. See for example \cite{G}.
A widely used  method  to discuss the local uniqueness
 is to prove the  uniqueness of solution for the reduced finite dimensional  problem by
counting the local degree. Such method involves the estimates of the second order derivatives of the solutions, which
are quite lengthy and technical. Let us point out that such
 method is hard to apply to \eqref{7-4-11}, because the solutions of \eqref{7-4-11} are not $C^2$ anymore.
In this paper, we will use the following
Pohozaev identities for the solution $u$ of $-\Delta u = f(x, u)$  to prove the local uniqueness result:
 \begin{equation}\label{1-1-11}
\begin{split}
& -\int_{\partial B_{\tau}(x_{0})}\frac{\partial u}{\partial \nu}\frac{\partial u}{\partial x_i}
+\frac12 \int_{\partial B_{ \tau}(x_{0})}|\nabla
u|^2 \nu_i\\
& =\int_{\partial B_{\tau}(x_{0})} F(x, u)\nu_i -\int_{ B_{\tau}(x_{0})} F_{x_i}(x, u),\quad i=1,\cdots, N,
\end{split}
\end{equation}
where  $\nu=(\nu_1,\cdots,\nu_N)$ is the outward unit normal of $\partial
B_{ \tau}(x_{0 })$,  $F(x, t)= \int_0^t f(x, s)\,ds$.
The advantage of such method is that we only need to estimate
the first order derivatives of the solutions, though this is not an easy task due to the jumping nonlinearities in  \eqref{7-4-11}.
The Pohozaev identities were used in \cite{DLY} to study the local uniqueness and periodicity of the solutions for the prescribed scalar equation.
Thanks to the coefficient in the prescribed scalar equation,  the Pohozaev identities \eqref{1-1-11} have a volume integral in the right hand side, which dominates all the surface integrals.
The estimate of such volume integral is relatively simple because it can be achieved by standard scaling argument. In the problem we consider now,
only line integrals appears in the Pohozaev identities. So we need to carefully study each line integral to determine which one dominates all the others.

\section{ Asymptotic of the  solutions}

Let  $\psi_{\lambda}$ be  a solution to \eqref{4-7-9} satisfying
\eqref{10-7-9}. We assume that as $\lambda\to +\infty$, $diam\,
\Omega_{\lambda, j}\to 0$. Throughout this section, we will denote $r_{\lambda, j} =\frac12 diam\,
\Omega_{\lambda, j}$ and let $p_{\lambda, j}\in \Omega_{\lambda, j} $ be
a point satisfying $\psi_{\lambda} (p_{{\lambda}, j})
=\max_{x\in \Omega_{\lambda, j}} \psi_{\lambda}(x)$.

\begin{lemma}\label{l1-28-7}

 Let  $\psi_{\lambda}$ be  a solution to \eqref{4-7-9}.  For  $x\in \Omega\setminus \cup_{j=1}^k \{ x:
d(x, \Omega_{\lambda, j}) \le L r_{\lambda, j}\}$,  it holds
\[
\psi_{\lambda}(x) =\sum_{j=1}^k \kappa_j G( p_{\lambda, j}, x) +
O\Bigl(\sum_{j=1}^k \frac{r_{\lambda, j}}{|x-p_{\lambda, j} |}
\Bigr),
\]
and
\[
\frac{\partial \psi_{\lambda}(x)}{\partial x_i} =\sum_{j=1}^k \kappa_j\frac{\partial G( p_{\lambda, j}, x)}{\partial x_i}  +
O\Bigl(\sum_{j=1}^k \frac{r_{\lambda, j}}{|x-p_{\lambda, j} |^2}
\Bigr),
\]
where $L>0$ is a large constant.

\end{lemma}

\begin{proof}

For any   $x\in  \Omega\setminus \cup_{j=1}^k \{ x:  d(x,
\Omega_{\lambda, j})\le L r_{\lambda, j}\}$, it holds $x\notin
\Omega_{\lambda, j}$.  Noting that
\[
H(y, x)- H( p_{\lambda, j}, x)=O\bigl(\frac{r_{\lambda, j}}{|x-p_{\lambda, j}|}\bigr), \quad y\in \Omega_{\lambda, j},
\]
we find
\[
\begin{split}
\psi_{\lambda}(x) =& \lambda \sum_{j=1}^k \int_{\Omega_{\lambda, j}} G(y, x)\,dy\\
=&\sum_{j=1}^k \lambda|\Omega_{\lambda, j}| G( p_{\lambda, j}, x)+
 \lambda \sum_{j=1}^k \int_{\Omega_{\lambda, j}} \Bigl( G(y, x)- G( p_{\lambda, j}, x)\Bigr)\,dy\\
=&\sum_{j=1}^k\kappa_j G( p_{\lambda, j}, x)+ \frac\lambda{2\pi}
\int_{\Omega_{\lambda, j}}\ln \frac{|x-p_{\lambda, j}|}{|y-x|}\,dy +
O\Bigl(\sum_{j=1}^k \frac{r_{\lambda, j}}{|x-p_{\lambda, j}|} \Bigr).
\end{split}
\]
Since
\begin{equation}\label{add2}
|y-x| = |x-p_{\lambda, j}|-\bigl\langle\frac{x-p_{\lambda, j}}{
|x-p_{\lambda, j}|}, y-p_{\lambda,j}\bigr\rangle+O\bigl(\frac{
|y-p_{\lambda, j}|^2}{ |x-p_{\lambda, j}|}\Bigr),\quad y\in
\Omega_{\lambda, j},
\end{equation}
 the result follows in $C(\Omega\setminus \cup_{j=1}^k \{ x:  d(x, \Omega_{\lambda, j}) \le L r_{\lambda, j}\})$.

Similarly, for $i=1,\cdots,k$,
 \[
\begin{split}
\frac{\partial \psi_{\lambda}(x)}{\partial x_i} =& \lambda \sum_{j=1}^k \int_{\Omega_{{\lambda}, j}} \frac{ \partial G(y, x)}{\partial x_i}\,dy\\
=&\sum_{j=1}^k\kappa_j \frac{\partial G( p_{{\lambda}, j}, x)}{\partial x_i}+
 \lambda \sum_{j=1}^k \int_{\Omega_{{\lambda}, j}} \Bigl(  \frac{\partial G(y, x)}{\partial x_i}- \frac{\partial G( p_{{\lambda}, j}, x)}{\partial x_i}\Bigr)\,dy\\
=&\sum_{j=1}^k\kappa_j\frac{\partial G(x, p_{{\lambda}, j})}{\partial x_i}+ \frac{\lambda}{2\pi} \sum_{j=1}^k\int_{\Omega_{{\lambda}, j}}\Bigl( \frac{x_i-p_{{\lambda},j,i}}{|x-p_{{\lambda},j}|^2}+\frac{y_i-x_i}{|y-x|^2} \Bigl)\,dy \\
& +  O\Bigl(\sum_{j=1}^k \frac{r_{\lambda, j}}{|x-p_{\lambda, j} |^2} \Bigr).
\end{split}
\]
Write
\[
\frac{x_i-p_{{\lambda},j,i}}{|x-p_{{\lambda},j}|^2}+\frac{y_i-x_i}{|y-x|^2}=\frac{y_i-p_{{\lambda},j,i}}{|x-p_{{\lambda},j}|^2}+ (y_i-x_i)\Bigl(\frac{1}{|y-x|^2}
-\frac{
1}{|x-p_{{\lambda},j}|^2}
\Bigr).
\]
Using again \eqref{add2}, we obtain
$$
\frac{\partial \psi_{\lambda}(x)}{\partial x_i}=\sum_{j=1}^k\kappa_j \frac{\partial G( p_{{\lambda}, j}, x)}{\partial x_i}+O\Bigl(\sum_{j=1}^k \frac{r_{{\lambda}, j} }{|x-p_{{\lambda}, j} |^2}   \Bigr),
$$
and thus complete our proof of Lemma~\ref{l1-28-7}.
\end{proof}

\begin{lemma}\label{l1-11-9}

Let  $\psi_\lambda$ be  a solution to \eqref{4-7-9} and satisfy
\eqref{5-7-9}. It holds
\[
\Omega_{\lambda, j} =B_\delta(x_{0, j})
 \cap\{\psi_\lambda>\tilde \kappa_{\lambda, j}\},
 \]
for some $\tilde \kappa_{\lambda, j}\to +\infty$ (as $\lambda\to +\infty$).

\end{lemma}

\begin{proof}
It follows from \eqref{5-7-9} that
\begin{equation}\label{6-7-9}
\lambda \psi_\lambda+ P+ \frac12 |\mathbf v|^2=c_{1,j},\quad \text{in}\;
\Omega_{\lambda, j},
\end{equation}
and
\begin{equation}\label{7-7-9}
P+ \frac12 |\mathbf v|^2=c_{2,j},
\end{equation}
for some constant $c_{2,j}$ in each connected component of
$B_\delta(x_{0, j})\setminus \Omega_{\lambda, j}$. Here $c_{2,j}$ may
depend on each connected component of $B_\delta(x_{0, j})\setminus
\Omega_{\lambda, j}$. By the continuity of $P+ \frac12 |\mathbf
v|^2$, we deduce from \eqref{6-7-9} and \eqref{7-7-9} that $\psi_\lambda$ is
a constant on each connected component of $\partial \Omega_{\lambda,
j}$. Let us also point out that $\Omega_{\lambda, j}$ must be simple
connected. If not, $\Omega_{\lambda, j}$ has an inner boundary
$\Gamma$, and $\Delta \psi_\lambda=0$ in $S$, which is the domain enclosed
by $\Gamma$. This will give $\psi_\lambda$ is a constant in $S$ and thus
$\nabla \psi_\lambda =0$ in $S$. We get a contradiction  by using the strong
maximum principle for the equation $-\Delta \psi_\lambda =\lambda $ in
$\Omega_{\lambda, j}$.  So $\partial \Omega_{\lambda, j}$ just has
one connected piece, on which $\psi_\lambda= \tilde \kappa_{\lambda,j}$ for some
constant $\tilde \kappa_{\lambda,j}>0$.  Using the maximum principle, we
conclude that $\psi_\lambda>\tilde \kappa_{\lambda,j}$ in $\Omega_{\lambda, j}$. We
claim that $ \tilde \kappa_{\lambda,j}\to +\infty$ as $\lambda\to+\infty$.
Firstly,  \eqref{3-7-9} is equivalent to   $\lambda
|\Omega_{\lambda, j}|=\kappa_j$. By  Lemma~\ref{l1-28-7}, for any
$M>0$ large,  $\psi_\lambda\ge M$ on $\partial B_\theta(x_{0, j})$ if
$\theta>0$ is small. By the maximum principle, $\psi_\lambda>M$ in
$B_\theta(x_{0, j})$.  But $\Omega_{\lambda, j}\subset
B_\theta(x_{0, j})$. It holds $\tilde \kappa_{\lambda,j}>M$.

We claim that $\psi_\lambda< \tilde \kappa_{\lambda,j}$ in $B_\delta(x_{0,
j})\setminus \Omega_{\lambda, j}$.  Indeed,  Lemma~\ref{l1-28-7}
implies
 $\psi_\lambda\le C$ on $\partial B_\delta(x_{0, j})$. Moreover, $\Delta \psi_\lambda=0$ in $ B_\delta(x_{0, j})\setminus \Omega_{\lambda, j}$. This
 gives $\psi_\lambda<\tilde \kappa_{\lambda,j}$ in $B_\delta(x_{0, j})\setminus \Omega_{\lambda, j}$.
So we have proved that $\Omega_{\lambda, j} =B_\delta(x_{0, j})
 \cap\{\psi_\lambda>\tilde \kappa_{\lambda,j}\}$.

\end{proof}

Next, we prove the following result.

\begin{proposition}\label{l2-11-9}

As $\lambda\to +\infty$,
 it holds  $r_{\lambda, i}^{-1} d(p_{\lambda, i}, \partial \Omega)\to +\infty  $,
 $i=1, \cdots, k$. Moreover,
 for $j\ne i$.
 \[
\frac{|p_{\lambda, i}-p_{\lambda, j}|}{\max (r_{\lambda, j},
\,r_{\lambda, i})}\to +\infty.
\]

\end{proposition}

To prove Proposition~\ref{l2-11-9}, we need to prove some
lemmas. To start with, we have

\begin{lemma}\label{l2-28-7}

It holds
\[
\lambda \sum_{i=1}^k
 \int_{\Omega_{\lambda , i}}
(\psi_\lambda -\tilde\kappa_{\lambda, i})_+=O(1).
\]

\end{lemma}

\begin{proof}

We have the following Pohozaev identity:
\begin{equation}\label{1-28-7}
\begin{split}
&\int_{\partial \Omega}\bigl\langle x-p_{\lambda , j}, \nabla
\psi_{\lambda} \bigr\rangle \frac{\partial \psi_\lambda }{\partial
\nu}-\frac12 \int_{\partial \Omega}\bigl\langle x-p_{\lambda, j},
\nu \bigr\rangle |\nabla
\psi_{\lambda}|^2\\
= &2 \lambda \sum_{i=1}^k
 \int_{\Omega_{\lambda , i}}
(\psi_\lambda -\tilde\kappa_{\lambda, i})_+.
\end{split}
\end{equation}

Using Lemma~\ref{l1-28-7}, we find that the left hand side of
\eqref{1-28-7} is bounded.  So the result follows.

\end{proof}

Now we study the local behaviors of $\psi_{\lambda}$ near
$p_{\lambda, j}$. Let $v_{\lambda} = \psi_{\lambda}
-\tilde\kappa_{\lambda, j}$. Then
\begin{equation}\label{2-28-7}
-\Delta v_{\lambda} = \lambda \Bigl( 1_{\{v_{\lambda}>0\}}   +\sum_{i\ne
j} 1_{B_\delta (x_{0, i}) }1_{ \{v_\lambda>\tilde
\kappa_{\lambda,i}-\tilde\kappa_{\lambda, j}\}}
 \Bigr),\quad \text{in}\; \Omega.
\end{equation}

Let $\Omega_j= \bigl\{ y:  \,r_{\lambda,j} y+p_{\lambda,
j}\in\Omega\bigr\}$ and $ f( x, v_\lambda)=\sum_{i\ne j} 1_{B_\delta (x_{0, i}) }1_{\{
v_\lambda>\tilde \kappa_{\lambda,i}-\tilde\kappa_{\lambda, j}\}} $.
Set $\bar v_{\lambda}(y)= v_{\lambda}( r_{\lambda,j} y+p_{\lambda,
j})$. We have
\begin{equation}\label{3-28-7}
-\Delta \bar v_{\lambda} =\lambda r^2_{\lambda, j} \Bigl(
1_{\{\bar v_{\lambda}>0\}}   +f( r_{\lambda,j} y+p_{\lambda, j},
\bar v_\lambda)
 \Bigr),\quad
\text{in}\; \Omega_j.
\end{equation}

 Let $w_{\lambda}
=\frac{1}{\lambda r_{\lambda, j}^2}\bar v_{\lambda} $. Then
\begin{equation}\label{4-28-7}
\begin{cases}
-\Delta w_{\lambda} =  1_{\{w_{\lambda}>0\}}+f( r_{\lambda,j}
y+p_{\lambda, j},
\lambda r_{\lambda, j}^2 w_\lambda),& \text{in}\; \Omega_j,\\
w_{\lambda}=-\frac{\tilde \kappa_{\lambda, j}}{\lambda r_{\lambda, j}^2}  , & \text{on}\; \partial \Omega_j.
\end{cases}
\end{equation}

\begin{lemma}\label{l3-28-7}

  For any $R>0$, there is a constant $ C>0$, depending on $R$, such that
\[
\|w_{\lambda}\|_{L^\infty(B_R(0)\cap\Omega_j)}\le C.
\]

\end{lemma}

\begin{proof}

First, we prove
\begin{equation}\label{1-29-7}
\int_{B_{r_{\lambda, j}^{-1} \delta}(0)\cap\Omega_j} (w_{\lambda})_+\le C.
\end{equation}
It follows from Lemma~\ref{l2-28-7} that
\[
\begin{split}
&2\lambda \int_{\Omega_{\lambda, j}} (\psi_{\lambda}-\tilde\kappa_{\lambda, j})_+\\
 =&2 \lambda r_{{\lambda}, j}^2 \int_{B_{r_{\lambda, j}^{-1} \delta}(0)\cap\Omega_j} (\tilde v_{\lambda})_+
 =2 \int_{B_{r_{{\lambda}, j}^{-1} \delta}(0)\cap\Omega_j} (w_{\lambda})_+.
 \end{split}
\]
Thus \eqref{1-29-7} follows.

By \eqref{1-29-7}, using the Morse iteration, we can prove
\begin{equation}\label{2-29-7}
\|(w_{\lambda})_+\|_{L^\infty(B_R(0)\cap\Omega_j)}\le C.
\end{equation}

Using the Harnack inequality, we can conclude
$\|w_{\lambda}\|_{L^\infty(B_R(0)\cap\Omega_j)}\le C$.  In fact,  we let $w_1 $
be a solution of
\begin{equation}\label{3-29-7}
\begin{cases}
-\Delta w_1 = 1_{\{w_{\lambda}>0\}}+f( r_{\lambda,j} y+p_{\lambda, j},
\lambda r_{\lambda, j}^2 w_\lambda) ,& \text{in}\; B_{R} (0)\cap\Omega_j,\\
w_1=0, & \text{on}\; \partial (B_{R} (0)\cap\Omega_j).
\end{cases}
\end{equation}
Then  $|w_1|\le C$.  Now $ w_2:= w_{\lambda} -w_1$ satisfies $\Delta
w_2=0$ and
\begin{equation}\label{5-29-7}
 \sup_{B_R(0)\cap\Omega_j} w_2\ge \sup_{B_R(0)\cap\Omega_j} w_{\lambda} -C \ge -C,
 \end{equation}
  since  $\sup_{B_R(0)\cap\Omega_j} w_{\lambda}
\ge 0$.

On the other hand,  by \eqref{2-29-7}, we have
\[
\sup_{B_R(0)\cap\Omega_j} w_2\le \sup_{B_R(0)\cap\Omega_j} w_{\lambda} +C < M,
\]
 for some large constant $M>0$. Thus,  $M- w_2$ is positive harmonic function. By Harnack inequality, there exists a constant $L>0$, such that
\[
 \sup_{B_R(0)\cap\Omega_j} (M-w_2 )\le L \inf_{B_R(0)\cap\Omega_j} (M-w_2),
 \]
 which, together with  \eqref{5-29-7},  gives
\[
 \inf_{B_R(0)\cap\Omega_j} w_2\ge  M - LM + L \sup_{B_R(0)\cap\Omega_j} w_2\ge - M'.
 \]

\end{proof}

\begin{lemma}\label{p1-3-8}

As $\lambda\to +\infty$, we have  $r_{\lambda, j }^{-1} d(p_{\lambda, j}, \partial \Omega)\to +\infty  $,
 $j=1, \cdots, k$,  and  $w_{\lambda} \to w$ in
$C^1_{loc}(\mathbb R^2)$, where
\[
w=
\begin{cases}
\ds\frac14(1-|x|^2), & |x|\le 1;\vspace{2mm}\\
\ds\frac12\ln\frac1{|x|}, & |x|\ge 1.
\end{cases}
\]
Moreover,
\[
\frac{|\Omega_{{\lambda}, j}|}{r_{\lambda, j}^2} \to \pi,
\]
and
\[
  \frac1{2\pi}\ln\frac1{r_{\lambda, j}}+   \sum_{i\ne j}
 G(p_{\lambda, j}, p_{\lambda, i})- \frac{\tilde\kappa_{\lambda, j}}{\pi \lambda
 r_{{\lambda}, j}^2} -H(p_{\lambda, j}, p_{\lambda, j})  \to  0.
 \]

\end{lemma}

\begin{proof}

For $x\in  \bigl(B_{\delta r_{\lambda, j}^{-1}}(0)\setminus B_L(0)\bigr)\cap \Omega_j $,
where $L>0$ is a large constant, it follows from Lemma~\ref{l1-28-7}
that
\[
\begin{split}
&w_{\lambda}=\frac{1}{\lambda r_{\lambda, j}^2}\Bigl( \psi_{\lambda}(r_{\lambda, j} x+p_{\lambda, j})-\tilde\kappa_{\lambda, j}
\Bigr)\\
 =& \frac{|\Omega_{\lambda, j}|}{r_{\lambda, j}^2}\Bigl(  \sum_{i=1}^k  G(r_{\lambda, j} x+p_{\lambda
 , j}, p_{\lambda, i})- \frac{\tilde\kappa_{\lambda, j}}{\lambda|\Omega_{\lambda, j}|}+  O\bigl( \frac1L\bigr)  \Bigr)
\\
=&\frac{|\Omega_{\lambda, j}|}{r_{\lambda, j}^2}\frac1{2\pi} \ln\frac1{|x|}\\
& +\frac{|\Omega_{\lambda, j}|}{r_{\lambda, j}^2}\Bigl(
\frac1{2\pi}\ln\frac1{r_{\lambda, j}}+  \sum_{i\ne j}
 G(r_{{\lambda}, j} x+p_{\lambda, j}, p_{\lambda, i})- \frac{\tilde\kappa_{\lambda,j}}{
 \lambda |\Omega_{\lambda, j}|} -H(r_{\lambda, j} x+p_{\lambda, j}, p_{\lambda, j})+  O\bigl( \frac1L\bigr)  \Bigr).
\end{split}
\]

 From $\frac{|\Omega_{\lambda, j}|}{r_{\lambda, j}^2}\le C$, we assume (up to a subsequence)
 that $\frac{|\Omega_{\lambda, j}|}{r_{\lambda, j}^2}\to t\in [0, +\infty)$.
 By Lemma~\ref{l3-28-7}, $|w_{\lambda}(x)|\le C$ for any $x\in B_R(0)\cap \Omega_j$, which implies
 \[
  \frac{|\Omega_{\lambda, j}|}{r_{\lambda, j}^2}\Bigl( \frac1{2\pi}\ln\frac1{r_{\lambda, j}}+   \sum_{i\ne j}
 G(p_{\lambda, j}, p_{\lambda, i})- \frac{\tilde\kappa_{\lambda, j}}{\lambda|\Omega_{\lambda, j}|} -H(p_{\lambda, j},
  p_{\lambda, j})  \Bigr)\to  \alpha_j\in (-\infty, +\infty).
 \]

 We have two possibilities:  (i)    $r_{\lambda, j }^{-1} d(p_{\lambda, j}, \partial \Omega)\to +\infty  $;  (ii)
  $r_{\lambda, j }^{-1} d(p_{\lambda, j}, \partial \Omega)\to a<+\infty $.
   We will prove that case (ii) can not occur.

Suppose that (i) occurs. Then from \eqref{4-28-7}, we have   $w_{\lambda} \to w$ in $C^1_{loc}(\mathbb
R^2)$, and $w$ satisfies
\begin{equation}\label{1-3-8}
\begin{cases}
-\Delta w =  1_{\{w>0\}}+\sum_{i=1}^m 1_{\{w>\beta_i\}},& \text{in}\; B_{R} (0),\\
w= \frac t{2\pi} \ln \frac1{|x|} +\alpha_j, & \text{in}\;
B_{R} (0)\setminus B_L(0),
\end{cases}
\end{equation}
where $R>>L>>1$ are two constants,  and $\beta_i=\frac1{\lambda
r_{\lambda, j}^2} \bigl(\tilde
\kappa_{\lambda,i}-\tilde\kappa_{\lambda, j}\bigr)\in [-\infty,
+\infty] $. Here, $1_{\{w_\lambda>-\infty\}}=1$ and
$1_{\{w_\lambda\ge\infty\}}=0$.

Since  $\Delta w\le 0$, $w$ attains its minimum at the boundary of $B_R(0)$. So
 $w(x)\ge \frac t{2\pi} \ln \frac1{R} +\alpha_j$  for all $x\in B_R(0)$.  Using the method of moving plane, we conclude that
 the solution of \eqref{1-3-8} must be radially symmetric, and thus $\{x: \; w>0\}$ is a disk.   Since
 $\min_{y\in \partial \Omega_{\lambda, j}}|y-p_{\lambda, j}|\le  r_{\lambda, j}$
and  $\max_{y\in \partial \Omega_{\lambda, j}}|y-p_{\lambda, j}|\ge
r_{\lambda, j}$,   we can find a $z_{\lambda,j}\in
\partial \Omega_{\lambda, j}$, such that
$|z_{\lambda,j}-p_{\lambda, j}|=  r_{\lambda, j}$. Let $y_{\lambda}
=\frac{z_{\lambda,j}-p_{\lambda, j}}{r_{\lambda, j}}$. Then
$|y_{\lambda}|=1$  and $r_{\lambda, j} y_{\lambda} +p_{\lambda, j}=
z_{\lambda,j}\in
\partial \Omega_{\lambda, j}$. Thus,
$w_{\lambda}(y_{\lambda})=0$. As a result, there  exists a $y$,
$|y|=1$, such that  $w(y)=0$.  So $\{x: \; w>0\}= B_1(0)$, which
gives
\[
w(x)= \frac14\bigl( 1- |x|^2\bigr),\quad x\in B_1(0),
\]
and $\Delta w= 0$ in $B_R(0)\setminus B_1(0)$.  Since $w\in C^1(B_R(0))$, we have $w(x) =\frac12 \ln\frac1{|x|}$  in $B_R(0)\setminus B_1(0)$.
Comparing this with \eqref{1-3-8}, we conclude
\[
\frac{|\Omega_{\lambda, j}|}{r_{\lambda, j}^2} \to t =\pi,
\]
and
\[
  \frac1{2\pi}\ln\frac1{r_{\lambda, j}}+   \sum_{i\ne j}
 G(p_{\lambda, j}, p_{\lambda, i})- \frac{\tilde\kappa_{\lambda, j}}{\lambda|\Omega_{{\lambda}, j}|} -H
 (p_{\lambda, j}, p_{\lambda, j})  \to  0.
 \]

 Suppose case (ii) occurs.  We first claim that
 as $\lambda\to +\infty$,
\begin{equation}\label{1-29-11}
\frac{\tilde \kappa_{\lambda, j}}{\lambda r_{\lambda, j}^2}\to +\infty.
\end{equation}

We argue by contradiction.  Suppose that as $\lambda\to +\infty$,
\begin{equation}\label{2-29-11}
\frac{\tilde \kappa_{\lambda, j}}{\lambda r_{\lambda, j}^2}\to a\in [0, +\infty).
\end{equation}
Similar to case (i), we find that
 $w_{\lambda} \to w$ in $C^1_{loc}(\mathbb
R^2_+)$ and after suitable translation and rotation, $w$ satisfies
\begin{equation}\label{10-29-11}
\begin{cases}
-\Delta w =  1_{\{w>0\}}+\sum_{i=1}^m 1_{\{w>\beta_i\}},& \text{in}\; B_{R} (0)\cap \mathbb R^2_+,\\
w= \frac t{2\pi} \ln \frac1{|x|} +\alpha_j, & \text{in}\; \bigl(
B_{R} (0)\setminus B_L(0)\bigr)\cap \mathbb R^2_+,\\
w(x_1, 0)=-a, & x_1\in (-R_1, R_2).
\end{cases}
\end{equation}
Comparing the last two relations in \eqref{10-29-11}, we find $t=0$ and $\alpha_j= -a$.  Thus, $w$ attains its minimum $-a$
in the whole region $\bigl(
B_{R} (0)\setminus B_L(0)\bigr)\cap \mathbb R^2_+$. This is a contradiction to the strong maximum principle. So we have proved \eqref{1-29-11}.

Let  $w_1$ be the solution of
\begin{equation}\label{1-18-10}
\begin{cases}
-\Delta w_{1} =  1_{\{w_{\lambda}>0\}}+f( r_{\lambda,j}
y+p_{\lambda, j},
\lambda r_{\lambda, j}^2 w_\lambda),& \text{in}\; \Omega_j,\\
w_{1}=1  , & \text{on}\; \partial \Omega_j.
\end{cases}
\end{equation}
Then, $w_1>0$.  Let $w_2 = w_{\lambda}-w_1$ satisfies  $\Delta w_2=0$ in  $\Omega_j$  and  $w_2=-\frac{\tilde \kappa_{\lambda, j}}{\lambda r_{\lambda, j}^2}-1$
on $\partial \Omega_j$.  So,  $w_2=-\frac{\tilde \kappa_{\lambda, j}}{\lambda r_{\lambda, j}^2}-1$
in $\ \Omega_j$. This gives
\begin{equation}\label{2-18-10}
w_{\lambda}= w_1-\frac{\tilde \kappa_{\lambda, j}}{\lambda r_{\lambda, j}^2}-1.
\end{equation}

Since  $r_{\lambda, j}^{-1} d(p_{\lambda, j}, \partial \Omega)\to C<+\infty  $  as $\lambda \to +\infty$,
 for $R>0$ large, it holds  $\partial \Omega_j \cap B_R(0)\ne \emptyset$.  Using the Harnack inequality, noting that $w_\lambda$ is bounded
in $ \Omega_j \cap B_R(0)$,  we deduce
\begin{equation}\label{3-18-10}
\sup_{ \Omega_i \cap B_R(0)} w_1 \le C \bigl(\inf_{ \Omega_i \cap B_R(0)} w_1+1\bigr)\le 2 C,
\end{equation}
which, together with \eqref{2-18-10}, gives
\begin{equation}\label{4-18-10}
w_{\lambda}\le C'-\frac{\tilde \kappa_{\lambda, j}}{\lambda r_{\lambda, j}^2}<0, \quad \text{in}\; \Omega_j \cap B_R(0).
\end{equation}
This is a contraction to  $w_\lambda(0)>0$.

\end{proof}

Now we are ready to prove Proposition~\ref{l2-11-9}.
\begin{proof}[Proof of Proposition~\ref{l2-11-9}]

We will argue by contradiction.  Fix $j$ and suppose that there are $j_1,
\cdots, j_k$ and $j_h\ne j$, $h=1, \cdots, k$,  such that
\[
\frac{|p_{\lambda, j_h}-p_{\lambda, j}|}{max (r_{\lambda, j},
\,r_{\lambda, j_h})}\le C.
\]
 for some $C>0$.
Without loss of generality, we assume that $r_{\lambda, j} \ge
\max_h r_{\lambda, j_h}$.  Otherwise, we will replace $j$ by some
$j_h$.  This implies $\Omega_{\lambda, j_h}\subset  B_{ R r_{\lambda, j}}(p_{\lambda,
 j})$  for some $j_h\ne j$.

It follows from Lemma~\ref{p1-3-8} that
\begin{equation}\label{1-12-9}
\frac{\partial \psi_\lambda}{\partial \nu }<0, \quad \text{in }\;
B_{ R r_{\lambda, j}}(p_{\lambda, j}),
\end{equation}
where $\nu=\frac{x-p_{\lambda, j}}{|x-p_{\lambda, j}|}$.

   On the other hand, by the maximum principle, it holds $\psi_\lambda\ge
 \kappa_{\lambda, j_h}$ in $\Omega_{\lambda, j_h}\subset  B_{ R r_{\lambda, j}}(p_{\lambda,
 j})$,  and  $\psi_\lambda\ge
 \kappa_{\lambda, j}$ in $\Omega_{\lambda, j}$. Noting $\Omega_{\lambda, j}\cap \Omega_{\lambda, j_h}=\emptyset$, we obtain   a contradiction to \eqref{1-12-9}.

 \end{proof}

By Lemma~\ref{p1-3-8}, we obtain the following local estimate for the solution $\psi_{\lambda}$:
\begin{equation}\label{nn1-4-8}
\psi_{\lambda}(x) = \lambda r_{\lambda, j}^2 \Bigl(
w\bigl(\frac{x-p_{\lambda, j}}{r_{\lambda, j}}\bigr)+o(1)\Bigr) +
\tilde\kappa_{\lambda,j}, \quad x\in B_{L r_{\lambda,
j}}(p_{\lambda, j}),
\end{equation}
and
\begin{equation}\label{1-4-8}
r_{\lambda, j}^2 \Bigl(\ln\frac1{r_{\lambda, j}}+ 2\pi\sum_{i\ne j}
 G(p_{\lambda, j}, p_{\lambda, i}) -2\pi H(p_{\lambda, j}, p_{\lambda, j})+o(1)\Bigr)=
 \frac{2\tilde\kappa_{\lambda, j}}
 {\lambda}.
 \end{equation}

Now we can calculate the local vorticity strength of the flow:
\begin{equation}\label{20-7-9}
\kappa_j=\lambda |\Omega_{\lambda, j}|= \frac{2(\pi+o(1))
\tilde\kappa_{\lambda, j} }{\ln\frac1{r_{\lambda, j}}+
2\pi\sum_{i\ne j}
 G(p_{\lambda, j}, p_{\lambda, i}) -2\pi H(p_{\lambda, j}, p_{\lambda, j})
+o(1)},
 \end{equation}
which implies
\begin{equation}\label{20-7-9'}
 \frac{4\pi\tilde\kappa_{\lambda, j}}{\ln \lambda}= \kappa_j +o(1).
 \end{equation}

We can also deduce from \eqref{1-4-8} that
\begin{equation}\label{n1-4-8}
r_{\lambda, j}= \frac{\sqrt{\kappa_{\lambda, j}}}{\sqrt{\pi \lambda}}\bigl( 1+ O\bigl(\frac1{\ln\lambda}\bigr)\bigr),
 \end{equation}
where $\kappa_{\lambda,j}=\frac{4\pi}{\ln\lambda}\tilde\kappa_{\lambda,j}$.

 Let
$$
\psi_\lambda=\frac{\ln \lambda}{4\pi}u_\lambda,\quad   \bar \lambda = \frac{4\pi \lambda}{\ln \lambda}.
$$
Then, $\kappa_{\lambda,j}\to \kappa_j$ as $\lambda\to +\infty$, and $u_\lambda$ satisfies
\begin{equation}\label{1'}
\begin{cases}
- \Delta u_\lambda=\bar\lambda\sum\limits_{j=1}^k 1_{B_\delta(x_{0, j})} 1_{\{ u_\lambda> \kappa_{\lambda,j}\}}, &\text{in}\;\Om,\\
u_\lambda=0, &\text{on}\; \partial\Om.
\end{cases}
\end{equation}
Moreover, it holds
\begin{equation}\label{2'}
  \lambda |\{u_\lambda>\kappa_{\lambda,j}\}|=\kappa_j.
\end{equation}

  From now on,   we will mainly investigate problem~\eqref{1'}.  Firstly,  we will discuss the global approximation for the solution of
\eqref{1'}

Let $R>0$ be a large constant, such that for any $x\in \Om$,
$\Om\subset B_R(x)$. Consider the following problem:
\begin{equation}\label{21}
\begin{cases}
-\Delta u=\bar \lambda1_{\{u>a\}},\;\; u>0, &\text{ in}\; B_R(0),\\
u=0, &\text{on}\;\partial B_R(0),
\end{cases}
\end{equation}
where
$a>0$ is a constant.
Then, \eqref{21} has a unique  solution $U_{\lambda,a}(y)$, which can be written as
\begin{equation}\label{1-5-1}
U_{\lambda,a}(y)=
\begin{cases}
a+\frac {\bar  \lambda}{4}\bigl(s_\lambda^2-|y|^2
\bigr), &  |y|\le s_\lambda,\vspace{0.2cm}\\
a\ln\frac {|y|} R/\ln \frac {s_\lambda}R, & s_\lambda\le |y|\le R,
\end{cases}
\end{equation}
where $s_\lambda$ is the constant, such that $U_{\lambda,a}\in C^1(B_R(0))$.
So, $s_\lambda$ satisfies
\begin{equation}\label{1-1-2}
-\frac{\bar \lambda s_\lambda }{2}=\frac{a}{s_\lambda\ln \frac {s_\lambda}R}.
\end{equation}
From
\begin{equation}\label{1-31-5}
s_\lambda\sqrt{\ln \frac R{s_\lambda}}=\sqrt{\frac{2a}{\bar \lambda}},
\end{equation}
we see  that  if $\lambda>0$ is large, \eqref{1-1-2} is uniquely solvable for $s_\lambda>0$ small.
Moreover, we have the following expansion for $s_\lambda$:
\begin{equation}\label{2-1-2}
s_\lambda=\frac{\sqrt{2a}}{\sqrt{\bar \lambda \ln\lambda}}\Bigl(1+O\bigl(\frac{\ln\ln\lambda}{\ln\lambda}\bigr)\Bigr)=
\frac{\sqrt{a}}{\sqrt{\pi \lambda }}\Bigl(1+O\bigl(\frac{\ln\ln\lambda}{\ln\lambda}\bigr)\Bigr).
\end{equation}

For any $x \in \Om$,  define
$U_{\lambda,x,a}(y)=U_{\lambda,a}(y-x)$. Because $U_{\lambda,x, a}(y)$
 does not satisfy the
zero boundary condition, we need to  make a  projection. Let
\begin{equation}\label{1-21-11}
PU_{\lambda,x,a}(y)= U_{\lambda,x,a}(y)-\frac a{\ln \frac{R}{s_\lambda}} g(y,x),
\end{equation}
where $g(y,x)$ satisfies
\[
\begin{cases}
- \Delta g=0,& \text{in } \; \Om,\vspace{0.2cm}\\
g=\ln\frac{R}{|y-x|}, &\text{on}\; \partial \Om.
\end{cases}
\]
It is easy to see that
\[
g(y,x)=\ln R +2\pi H(y,x),
\]
where $H(y,x)$ is the regular part of the Green function.

 For each local maximum point $p_{\lambda, j}$, we choose  $x_{\lambda, j}\in B_\delta(p_{\lambda, j})$, which is to be determined later.
Let  ${\bf x}_{\lambda}=( x_{\lambda, 1}, \cdots, x_{\lambda, k})$. For ${\bf a}=(a_1,\cdots,a_k)$
denote
\begin{equation}\label{1-5-11}
\mathcal U_{\lambda, {\bf x}_\lambda, {\bf a}}= \sum_{j=1}^k PU_{\lambda, x_{\lambda,  j},a_{ j}}.
\end{equation}
where  $a_j$ is chosen suitably close to $\kappa_{j}$.

We will  choose ${\bf x}_\lambda$,  ${\bf a}_\lambda=(a_{\lambda,1},\cdots,a_{\lambda, k})$   and $s_{\lambda, j}$ , such that the following conditions hold:
\begin{equation}\label{1-16-8}
\nabla \mathcal U_{\lambda, {\bf x}_\lambda, {\bf a}_\lambda}( p_{\lambda, j} )=0,
\end{equation}
\begin{equation}\label{2-17-10}
a_{\lambda,i}= \kappa_{\lambda,i} + \frac{a_{\lambda,i}}{\ln\frac R {s_{\lambda, i}} }  g(x_{\lambda, i},x_{\lambda, i}) -\sum_{j\ne i}
\frac{a_{\lambda,j}}{\ln\frac R{s_{ \lambda, j}} }  \bar G(x_{\lambda, i},x_{\lambda, j}),
\end{equation}
and
\begin{equation}\label{1-6-11}
s_{\lambda,i}\sqrt{\ln \frac R{s_{\lambda, i}}}=\sqrt{\frac{2a_{\lambda, i}}{ \bar\lambda}},
\end{equation}
where $\bar G(y,x)=\ln \frac R{|y-x|} -g(y,x)$.  Note that $\bar G(y,x)= 2\pi G(y, x)$  and $G(y, x)$ is
the Green function of $-\Delta$ subject to the zero boundary condition.

Note that \eqref{1-16-8} can be written as
\begin{equation}\label{2-16-8}
- \frac{\bar  \lambda(p_{\lambda,i}-x_{\lambda, i} ) }{2}= \frac{a_{\lambda,i}}{\ln\frac R {s_{\lambda, i}} }\nabla  g(x_{\lambda, i},x_{\lambda, i}) -\sum_{j\ne i}
\frac{a_{\lambda,j}}{\ln\frac R{s_{ \lambda, j}} } \nabla \bar G(x_{\lambda, i},x_{\lambda, j}).
\end{equation}
We can solve \eqref{1-16-8},
\eqref{2-17-10} and \eqref{1-6-11} to obtain  $x_{\lambda, i}$, $s_{\lambda, i}$  and  $a_{\lambda,i}$,  $i=1, \cdots, k$.
Moreover, we have
\begin{equation}\label{3-16-8}
 |x_{\lambda,i}-p_{\lambda, i}| = O\bigl( \frac{1}{\lambda }  \bigr),
\end{equation}
\begin{equation}\label{1-20-10}
a_{\lambda,i}=\frac{\kappa_{\lambda,i} - \sum_{j\ne i}\frac{1}{\ln\frac R{s_{\lambda, j}} }
a_{\lambda,j}\bar G(x_{\lambda, i},x_{\lambda, j}) }{1-\frac{g(x_{\lambda, i},x_{\lambda, i})}{\ln\frac
R{s_{\lambda, i}}}}= \kappa_{\lambda,i} +O\bigl(\frac1{\ln \lambda}\bigr),
\end{equation}
and by \eqref{2-1-2} and \eqref{n1-4-8}, we find
\begin{equation}\label{5-16-8}
 |r_{\lambda,i}-s_{\lambda, i}| = O\Bigl( \frac{\ln\ln\lambda}{\lambda\ln\lambda}\Bigr).
\end{equation}

We will estimate
\begin{equation}\label{3-1-2}
\omega_\lambda= u_\lambda -\mathcal U_{\lambda, {\bf x}_\lambda, {\bf a}_\lambda}.
\end{equation}

\begin{lemma}\label{l1-16-8}

 As $\lambda \to +\infty$,
\[
\|\omega_\lambda\|_{L^\infty(\Omega)} \to 0.
\]

\end{lemma}

\begin{proof}

 Using \eqref{nn1-4-8}, \eqref{3-16-8}, \eqref{1-20-10}  and \eqref{5-16-8}, we can easily prove
\begin{equation}\label{6-16-8}
 \omega_\lambda \to 0,
\end{equation}
uniformly in $\cup_{j=1}^kB_{L r_{\lambda, j}}(p_{\lambda, j})$.

On the other hand, noting that  $H(x, p_{\lambda, j})\ge -C$  in $\Omega\cap B_\delta (p_{\lambda, j})\setminus  B_{L r_{\lambda, j}}(p_{\lambda, j})$,
 we find from Lemma~\ref{l1-28-7} and \eqref{20-7-9} that for $x\in \Omega\cap B_\delta (p_{\lambda, j})\setminus  B_{L r_{\lambda, j}}(p_{\lambda, j})$,
\begin{equation}\label{1-17-8}
\begin{split}
\psi_\lambda(x)\le &  \kappa_j \Bigl( \frac1{2\pi} \ln \frac1{|x- p_{\lambda, j}|} +O(1)\Bigr)
 +O(1)\\
 =& \frac{2\pi
\tilde\kappa_{\lambda, j} }{\ln\frac1{r_{\lambda, j}}}\Bigl( 1+ O\bigl(\frac1{ \ln\frac1{r_{\lambda, j}}   }\bigr)\Bigr)\Bigl( \frac1{2\pi} \ln \frac1{|x- p_{\lambda, j}|} +O(1)\Bigr)
 +O(1)\\
 \le &\tilde \kappa_{\lambda,j} \Bigl( 1-\frac{ \ln L }{|\ln r_{\lambda, j}|  } + \frac{ C }{|\ln r_{\lambda, j}|  }\Bigr)+O(1)<\tilde \kappa_{\lambda,j},
 \end{split}
 \end{equation}
if $L>0$ is large,  which gives
 \begin{equation}\label{4-16-8}
 u_\lambda(x)=\frac{4\pi \psi_\lambda}{\ln \lambda}<\kappa_{\lambda,j},\quad x\in \Omega\cap B_\delta (p_{\lambda, j})\setminus  B_{L r_{\lambda, j}}(p_{\lambda, j}).
\end{equation}
   Similar to \eqref{15-2-5}  and \eqref{10-29-5} in the Appendix,  we can  show that
\begin{equation}\label{4-17-8}
\mathcal U_{\lambda, {\bf x}_\lambda, {\bf a}_\lambda}(x) <\kappa_{\lambda,j}, \quad x \in \Omega\cap B_\delta (p_{\lambda, j})\setminus  B_{L r_{\lambda, j}}(p_{\lambda, j}).
 \end{equation}
As a result,
\begin{equation}\label{5-17-8}
\Delta \omega_\lambda=0, \quad \text{in}\; \Omega\setminus \cup_{j=1}^k    B_{L r_{\ep, j}}(p_{\lambda, j}).
 \end{equation}
By the maximum principle, it holds
\begin{equation}\label{6-17-8}
\|\omega_\lambda\|_{L^\infty(  \Omega\setminus \cup_{j=1}^k    B_{L r_{\lambda, j}}(p_{\lambda, j})  )}\le
\|\omega_\lambda\|_{L^\infty(  \cup_{j=1}^k   \partial B_{L r_{\lambda, j}}(p_{\lambda, j})  )} \to 0.
 \end{equation}

\end{proof}

\section{the estimate of the error term}

Let
\begin{equation}\label{30}
w(y)=
\begin{cases}
\ds\frac14(1-|y|^2), &|y|\le 1,\vspace{0.2cm}\\
\ds\frac12\ln \frac1{|y|},  & |y|>1.
\end{cases}
\end{equation}
Then $w\in C^1(\mathbb R^2)$. It is easy to check  that $w$ satisfies
\begin{equation}\label{31}
-\Delta w= 1_{\{w>0\}} \quad \text{in}\; \mathbb R^2.
\end{equation}
Note that $w>0$ if $|y|<1$ and $w<0$ if $|y|>1$.

The linearized  operator for \eqref{31} is
\begin{equation}\label{32}
-\Delta v-  2v(1,\theta)\delta_{|y|=1}=0.
\end{equation}

We have  proved in \cite{CPY1} the following result:

\begin{proposition}\label{p31}

Let $v\in L^\infty(\mathbb R^2)\cap  C(\mathbb R^2) $ be a solution of \eqref{32}. Then
\[
v  \in span\bigl\{\frac{\partial w}{\partial y_1}, \;\frac{\partial
w}{\partial y_2}\bigr\}.
\]

\end{proposition}

Define the linear operator $\mathbb L_\lambda$ as follows.
\begin{equation}\label{10-17-8}
\mathbb L_\lambda  \omega= -  \Delta \omega- 2\sum_{j=1}^k \frac{1}{s_{\lambda, j} }
\omega(s_{\lambda, j},\theta)\delta_{|y-x_{\lambda, j}|=s_{\lambda, j}}, \quad \omega\in W^{1, p}_0(\Omega),
\end{equation}
where $p>2$.

Note that for any $\omega\in C(\Omega)$,
$
\omega(s_{\lambda, j},\theta)\delta_{|y-x_{\lambda, j}|=s_{\lambda, j}}
$
is a  bounded linear functional in $W^{1,q}(\Omega)$  for any $q\ge 1$.
  We have

\begin{proposition}\label{p32}
  For any  $p\in (2, +\infty]$,
 there are constants $c_0>0$ and $\lambda_0>0$, such that for any
$\lambda\ge\lambda_0$,
  $ \omega $ with  $\nabla \omega(p_{\lambda, j})=0$,  $j=1,\cdots,k$,  and
$ \mathbb L_\lambda \omega =0$ in $\Omega \setminus \cup_{j=1}^k
 B_{L s_{\lambda, j}}
(x_{\lambda, j})$  for some large $L>0$,  there holds
\[
\begin{split}
&s_{\lambda, 1}^{\frac2{p'}-1 }\| \mathbb L_\lambda \omega\|_{W^{-1, p}(\cup_{j=1}^k
 B_{L s_{\lambda, j}}
(x_{\lambda, j}))}+\| \mathbb L_\lambda \omega\|_{L^\infty(\cup_{j=1}^k
 B_{\frac12 s_{\lambda, j}}
(x_{\lambda, j}))}
 \\
 \ge & c_0 \Bigl( s_{\lambda, 1}^{1-\frac2p }
\|\nabla \omega\|_{L^{ p}(\cup_{j=1}^k
 B_{L s_{\lambda, j}}
(x_{\lambda, j}))} +\|\omega\|_{L^\infty(\Om)}\Bigr).
\end{split}
\]

\end{proposition}

\begin{proof}

We argue by contradiction. Suppose that there are $\{\lambda_n\}$,
$p_{n, j}\in \Omega$, and $\omega_n\in W^{1, p}_0(\Omega)$ with $\lambda_n\to +\infty$,  $\nabla\omega_n  (p_{n, j})=0$ and
$\mathbb L_{\lambda_n} \omega_n =0$ in $\Omega \setminus
  \cup_{j=1}^k
 B_{L s_{\lambda_n, j}}
(x_{n, j})  $, such that
\begin{equation}\label{1-8-11}
s_{\lambda_n, 1}^{\frac2{p'}-1} \| \mathbb L_{\lambda_n} \omega_n\|_{W^{-1, p}(\cup_{j=1}^k
 B_{L s_{\lambda_n, j}}
(x_{n, j}))} +\| \mathbb L_{\lambda_n} \omega\|_{L^\infty(\cup_{j=1}^k
 B_{\frac12 s_{\lambda_n, j}}
(x_{\lambda, j}))}\le
\frac{1}{n},
\end{equation}
and
\begin{equation}\label{2-11-11}
 s_{\lambda_n, 1}^{1-\frac2p }\|\nabla \omega_n\|_{L^{
p}(\cup_{j=1}^k
 B_{L s_{\lambda_n, j}}
(x_{n, j}))}+\|\omega_n\|_{L^\infty(\Omega)} =1.
\end{equation}

Let  $f_n= \mathbb L_{\lambda_n} \omega_n$,  then
\begin{equation}\label{21-24-5}
-\Delta \omega_n =\frac{2}{s_{\lambda_n, j} }
\omega_n(s_{\lambda_n, j},\theta)\delta_{|y-x_{\lambda_n, j}|=s_{\lambda_n, j}}+f_n.
\end{equation}

  Let  $\varphi_{n,j}(y)= \omega_n( s_{\lambda_n, j} y+ x_{\lambda_n, j})    $,  then
\begin{equation}\label{100-7-11}
\begin{split}
\int_{\mathbb R^2} \nabla
 \varphi_{n,j}\nabla \phi= 2\int_{|y|=1}
\varphi_{n,j}\phi +\bigl\langle  f_{n}, \phi\bigl(\frac{ y- x_{\lambda_n, j}  }{s_{\lambda_n, j}   }\bigr)\bigr\rangle,\quad \forall\; \phi\in C^\infty_0(\mathbb R^2).
\end{split}
\end{equation}

Since  the right hand side of \eqref{100-7-11} is bounded in $W^{-1,p}_{loc}
(\mathbb R^2)$ , $\varphi_{n,j}$ is bounded in $W^{1,p}_{loc}
(\mathbb R^2)$.
Noting that $p>2$, we deduce  from the
Sobolev embedding that $\varphi_{n,j}$ is bounded in $C^\alpha_{loc}(\mathbb R^2)$
for some $\alpha>0$. So, we can
assume that $\varphi_{n,j}$ converges uniformly in any compact set of $\mathbb R^2$
 to $\omega\in L^\infty(\mathbb R^2)\cap C(\mathbb R^2)$. It is
easy to check that  $\omega$ satisfies
\begin{equation}\label{22-24-5}
-\Delta \omega = 2\omega (1,\theta)\delta_{|y|=1}, \quad \text{in}\; \mathbb R^2.
\end{equation}
So, by Proposition~\ref{p31},
\begin{equation}\label{22-24-5}
\omega= c_1 \frac{\partial w}{\partial x_1} + c_2 \frac{\partial w}{\partial x_2}.
\end{equation}

On the other hand,  from $|f_n|\le \frac1n$ in $B_{\frac12 s_{\lambda, j}}
(x_{\lambda, j})$  and  $|\varphi_{n,j}|\le 1$,  we can deduce $\varphi_{n,j}$  is bounded in $W^{2, p}(B_{\frac14}(0))$. So
we can also assume $\varphi_{n,j}\to \omega$  in $ C^1(B_{\frac14}(0))$.
  Since $\nabla \varphi_{n,j} \bigl(  \frac{p_{\lambda, j}-x_{\lambda, j}}{ s_{\lambda, j}}  \bigr)=s_{\lambda, j}\nabla \omega_n(p_{\lambda, j})=0$
  and  $\frac{p_{\lambda, j}-x_{\lambda, j}}{ s_{\lambda, j}} \to 0$, it holds $\nabla\omega  (0)=0$. This implies
  $c_1= c_2=0$. That is, $\omega \equiv
0$. Thus, we have proved
\[
\omega _n =o(1),\quad \text{in}\; B_{L s_{\lambda_n, j}}(x_{n, j}),
\]
for any $L>0$.

By our assumption,
\[
\mathbb L_{\lambda_n} \omega_n = 0,\quad\text{in}\; \Om\setminus
\cup_{j=1}^k B_{L s_{\lambda_n, j}}(x_{n, j}).
\]
Thus,
\[
\Delta \omega_n  =0,\quad y\in \Om\setminus
\cup_{j=1}^k B_{L s_{\lambda_n, j}}(x_{n, j}).
\]
However, $\omega_n=0$ on $\partial\Om$      and $\omega_n=o(1)$ on
$\partial B_{L s_{\lambda_n, j}}(x_{n, j})$. By the maximum principle,
\[
\omega_n=o(1).
\]
So, we have proved that
\begin{equation}\label{327}
\|\omega_n\|_{L^\infty(\Omega)}=o(1)\,\,\,\hbox {as}\,\,\,n\to +\infty.
\end{equation}
 Moreover, it follows from \eqref{100-7-11} and the
Sobolev embedding that for any $\phi\in C_0(B_{2L }(0))$,
\begin{equation}\label{1-28-5}
\begin{split}
&\Bigl| \int_{\mathbb R^2}\nabla \varphi_{n,j} \nabla\phi\Bigr|\\
=&
\Bigl|2\int_{|y|=1} \varphi_{n,j}(1,\theta) \phi(1,\theta)+\bigl\langle  f_{n}, \phi\bigl(\frac{ y- x_{\lambda_n, j}  }{s_{\lambda_n, j}   }\bigr)\bigr\rangle
 \Bigr|\\
=& o(1)\int_{|y|=1}| \phi(1,\theta)| +o(1)  \|\phi\|_{W^{1, p'}(
B_{2L}(0))}\\
=& o(1)\|\phi\|_{W^{1, 1}(B_1(0))}  +o(1)\|\phi\|_{W^{1, p'}(
B_{2L}(0))}\\
=& o(1)\Bigl(\int_{B_{2L}(0)}
|\nabla\phi|^{p'}\Bigr)^{\frac1{p'}},
\end{split}
\end{equation}
which implies
\[
\|\nabla \varphi_{n,j} \|_{L^p(B_{2L}(0))}=o(1).
\]
That is,
\begin{equation}\label{1-11-11}
s_{\lambda_n, j}^{1-\frac2{p}}\|\nabla  \omega_{n} \|_{L^p(B_{2Ls_{\lambda_n,
j}}(x_{n, j}))}=o(1), \quad j=1, \cdots, k.
\end{equation}

Noting that  $0<c_0 \le \frac{s_{\ep_n, 1}}{s_{\ep_n, j}}\le c_1<+\infty$, we obtain a contraction from
\eqref{2-11-11}, \eqref{327} and \eqref{1-11-11} and thus complete our proof of Proposition~\ref{l2-11-9}.

\end{proof}

We are now ready to estimate the error term $\omega_\lambda$  defined by \eqref{3-1-2}.
  Note that  $\omega_\lambda$ satisfies
\begin{equation}\label{15-10-9}
 \mathbb L_\lambda \omega_\lambda=
   R_\lambda(\omega_\lambda),
\end{equation}
where
\begin{equation}\label{6-10-9}
\begin{split}
R_\lambda(\omega_\lambda)=&\bar \lambda\Bigl( \sum_{j=1}^ k  1_{B_\delta(x_{0, j})}  1_{\{ \mathcal U_{\lambda, {\bf x}_\lambda, { \bf a}_\lambda} +\omega_\lambda>\kappa_{
\lambda,j}\}}
-  \sum_{j=1}^ k 1_{\{U_{\lambda, x_{\lambda,j}, a_{\lambda, j}}>a_{\lambda, j} \} }\Bigr)\\
&-2\sum_{j=1}^k \frac{1}{s_{\lambda, j} }
\omega_\lambda(s_{\lambda, j},\theta)\delta_{|y-x_{\lambda, j}|=s_{\lambda, j}}.
\end{split}
\end{equation}
Then, from \eqref{4-16-8}  and \eqref{4-17-8}, it holds
\begin{equation}\label{1-6-2}
 R_\lambda(\omega_\lambda)=0\quad \text{in}\; \Omega \setminus \cup_{j=1}^k B_{2Ls_{\lambda,j}}(x_{\lambda, j}).
\end{equation}

By the choice of $x_{\lambda, j}$, we have
\begin{equation}\label{20-17-8}
 \nabla \omega_\lambda(x_{\lambda, j})=0.
\end{equation}

\begin{proposition}\label{p34}
 Fix a constant $p>2$.
Then
\[
s_{\lambda, j}^{1-\frac2p}\|\nabla
\omega_{\lambda}\|_{L^{p}(\cup_{j=1}^k B_{2Ls_{\lambda,j}}(x_{\lambda, j}))}  + \|
\omega_{\lambda}\|_{L^\infty(\Omega)}
= O\bigl(\sum_{j=1}^k \frac{s_{\lambda, j}}{
|\ln s_{\lambda, j}|}\bigr).
\]

\end{proposition}

\begin{proof}

 By \eqref{1-6-2} and \eqref{20-17-8},  we can apply Proposition~\ref{p32} to obtain
\begin{equation}\label{35-29-5}
\begin{split}
&s_{\lambda, 1}^{1-\frac2p}\|\nabla
\omega_\lambda\|_{L^{p}(\cup_{j=1}^k B_{2Ls_{\lambda,j}}(x_{\lambda, j}))}  + \|
\omega_\lambda\|_{L^\infty(\Omega)}\\
\le  & C s_{\lambda, 1}^{\frac2{p'}-1}\|
   R_\lambda(\omega_\lambda)\|_{W^{-1,p}(\cup_{j=1}^k B_{2Ls_{\lambda,j}}(x_{\lambda, j}))
  }+C
 \|
   R_\lambda(\omega_\lambda)\|_{L^\infty(\cup_{j=1}^k B_{\frac12s_{\lambda,j}}(x_{\lambda, j}))
  } \\
  =& C s_{\lambda, 1}^{\frac2{p'}-1}\|
   R_\lambda(\omega_\lambda)\|_{W^{-1,p}(\cup_{j=1}^k B_{2Ls_{\lambda,j}}(x_{\lambda, j}))
  }
  \end{split}
\end{equation}
since  Lemmas~\ref{l1-23-5}  and \ref{al1} gives
\begin{equation}\label{10-20-10}
   R_\lambda(\omega_\lambda)=0,  \quad \text{in}\; \cup_{j=1}^k B_{\frac12s_{\lambda,j}}(x_{\lambda, j}).
\end{equation}

Now, we estimate $s_{\lambda, 1}^{\frac2{p'}-1} \| R_\lambda(\omega)\|_{W^{-1,p}(\cup_{j=1}^k B_{2Ls_{\lambda,j}}(x_{\lambda, j}))} $.

  For  $j$, we use  $\tilde  \xi_j(y)$ to denote $\xi(s_{\lambda, j} y+ x_{\lambda, j})$  for any function  $\xi$.   Let $\phi\in C^1_0(B_{2Ls_{\lambda,j}}(x_{\lambda, j}
  ))$.  Then,
\begin{equation}\label{37-29-5}
\bigl\langle  R_\lambda(\omega_\lambda), \phi\bigr\rangle=
 s_{\lambda, j}^2\bar \lambda\int_{B_{2L}(0)}\Bigl(  1_{\{\tilde  {\mathcal U}_{\lambda, {\bf x}_\lambda, {\bf a}_\lambda,j}
 +\tilde \omega_{\lambda,j}>\kappa_{\lambda,j}\}} -1_{\{\tilde U_{\lambda,x_{\lambda,j},a_{\lambda, j},j}> a_{\lambda, j}\}}
  \Bigr)\tilde\phi_j -2 \int_{|y|=1}\tilde \omega_{\lambda,j}\tilde  \phi_j.
\end{equation}

Denote  $y_{\lambda, j}(\theta) =( (1+\tilde t_{\lambda, j}(\theta) ) \cos \theta, (1+\tilde t_{\lambda, j}(\theta) ) \sin\theta)$,
where $\tilde t_{\lambda, j}(\theta)$ is defined in Lemma~\ref{al1}.
Then, following Lemma~\ref{al1}, we find
\begin{equation}\label{1-18-11}
| |y_{\lambda, j}(\theta)|-1|= O\Bigl( |\tilde \omega_{\lambda, j}(y_{\lambda, j}(\theta)||\ln s_{\lambda, j}| +\frac1\lambda\Bigr).
\end{equation}

 It  follows from Lemma~\ref{al1} that
\begin{equation}\label{n39-29-5}
\begin{split}
&s_{\lambda, j}^2\bar \lambda\int_{B_{2L}(0)}\Bigl(  1_{\{\tilde  {\mathcal U}_{\lambda, {\bf x}_\lambda, {\bf a}_\lambda, j}
 +\tilde \omega_{\lambda, j}>\kappa_{\lambda,j}\}} -1_{\{\tilde U_{\lambda,x_{\lambda,j},a_{\lambda, j}, j}> a_{\lambda, j}\}}
  \Bigr)\tilde\phi_j
\\
=&2\pi s_{\lambda, j}^2\bar \lambda
\int_{0}^{2\pi}\int_1^{1+\tilde t_{\lambda, j}(\theta)}
\tilde \phi_j( t, \theta)\,dt
d\theta
\\
=& 2\pi s_{\lambda, j}^2\bar \lambda\int_{0}^{2\pi} \int_1^{1+\tilde t_{\lambda, j}(\theta)}
\tilde \phi_j( 1, \theta)\,dt
d\theta\\
&+2\pi s_{\lambda, j}^2\bar \lambda\int_{0}^{2\pi}\int_1^{1+\tilde t_{\lambda, j}(\theta)} \Bigl(
\tilde \phi_j(t,\theta)-\tilde \phi_j( 1, \theta)\Bigr)\,dt
d\theta\\
=&2\pi s_{\lambda, j}^2\bar \lambda\int_{0}^{2\pi} \Bigl(\frac{\tilde \omega_{\lambda,j}(y_{\lambda, j}(\theta) )}{a_{\lambda, j}}\ln\frac R{s_{\lambda, j}}- s_{\lambda, j}  \bigl\langle
\nabla g(x_{\lambda, j},x_{\lambda, j}\\
&\qquad\qquad  -\sum_{i\ne j}\frac{ a_{\lambda, j}}{a_{\lambda, j}}\nabla \bar G(x_{\lambda, i},x_{\lambda, j}) , (\cos\theta, \sin\theta)\bigr\rangle  \Bigr)\tilde\phi_j( 1, \theta)
\,d\theta\\
&+ O\Bigl(\frac1\lambda+ |\tilde \omega_{\lambda, j} (y_{\lambda, j})|^2 |\ln s_{\lambda, j}|\Bigr)\int_{|y|=1} |\tilde \phi_j|
\\
&+2\pi s_{\lambda, j}^2\bar \lambda\int_{0}^{2\pi}\int_1^{1+\tilde t_{\lambda, j}(\theta)}\Bigl(
\tilde \phi_j(t,\theta)-\tilde \phi_j( 1, \theta)\Bigr)\,dt
d\theta\\
=&2  \int_{|y|=1}\tilde  \omega_{\lambda, j} \tilde \phi_j+  4\pi   \int_{0}^{2\pi}\bigl(\tilde  \omega_{\lambda, j}(y_{\lambda, j}(\theta))- \tilde  \omega_{\lambda, j}( 1, \theta)\bigr) \tilde \phi_j( 1, \theta)\,d\theta
\\
&+\Bigl[ O\Bigl(\sum_{i=1}^k\frac{s_{\lambda, i}}{|\ln s_{\lambda, i}|}\Bigr)+ o(1)\|\omega_\lambda\|_{L^\infty(\Omega)}\Bigr]
\|\tilde \phi_j\|_{W^{1,p'}(B_{2L}(0))}\\
&+2\pi s_{\lambda, j}^2\bar \lambda\int_{0}^{2\pi}\int_1^{1+\tilde t_{\lambda, j}(\theta)}\Bigl(
\tilde \phi_j(t,\theta)-\tilde \phi_j( 1, \theta)\Bigr)\,dt
d\theta.
\end{split}
\end{equation}
Moreover, from \eqref{12-29-5}  and \eqref{1-18-11},
\begin{equation}\label{nn39-29-5}
\begin{split}
&2\pi s_{\lambda, j}^2\bar \lambda\int_{0}^{2\pi}\int_1^{1+\tilde t_{\lambda, j}(\theta)}\Bigl(
\tilde \phi_j(t,\theta)-\tilde \phi_j( 1, \theta)\Bigr)\,dt
d\theta\\
=&2\pi  s_{\lambda, j}^2\bar \lambda\int_{0}^{2\pi}\int_1^{1+\tilde t_{\lambda, j}(\theta)}
\int_1^ t
\frac{\partial  \tilde\phi_j(s,\theta)}{\partial s}\,ds dt d\theta\\
=& O\Bigl( s_{\lambda, j}^2\bar \lambda\int_{0}^{2\pi}|\tilde t_{\lambda, j}(\theta)|
\int_1^ {1+
|\tilde t_{\lambda, j}(\theta)|}
\bigl|\frac{\partial  \tilde\phi_j(s,\theta)}{\partial s}\bigr|\,ds  d\theta
\Bigr)\\
=&  O\Bigl(s_{\lambda, j}^2\bar \lambda( |\tilde \omega_{\lambda, j}(y_{\lambda, j}(\theta))||\ln s_{\lambda, j}| +s_{\lambda, j})\Bigr)
\int_{0}^{2\pi}
\int_1^ {1+
|\tilde t_{\lambda, j}(\theta)|}
\bigl|\frac{\partial  \tilde\phi_j(s,\theta)}{\partial s}\bigr|\,ds  d\theta\Bigr)\\
=&\Bigl[O \bigl(  s_{\lambda, j}^{1+\frac1p}+o(1)\| \omega_{ \lambda}\|_{L^\infty(\Omega)}
\Bigr]\| \tilde\phi_j\|_{W^{1,p'}(B_{2L}(0))}.
\end{split}
\end{equation}
We also have
\begin{equation}\label{30-17-8}
\begin{split}
&   \int_{0}^{2\pi}\bigl(\tilde  \omega_{\lambda, j}(y_{\lambda, j}(\theta))- \tilde  \omega_{\lambda, j}(1, \theta)\bigr) \tilde \phi_j(1, \theta)
\\
=& \int_{0}^{2\pi}\int_1^{1+\tilde t_{\lambda, j}(\theta)}
\frac{\partial  \tilde\omega_{\lambda, j}(s,\theta)}{\partial s}\tilde \phi_j(1, \theta)\,ds  d\theta\\
=& O\Bigl( \|\nabla \tilde \omega_{\lambda, j}\|_{L^{ p}(B_L(0))}\Bigl(\int_{0}^{2\pi}
\int_1^ {1+
|\tilde t_{\lambda, j}(\theta)|}
|\tilde \phi_j(1, \theta)|^{p'}\,ds  d\theta\Bigr)^{\frac1{p'}}
\Bigr)\\
=&o(1)\|\nabla \tilde \omega_{\lambda, j}\|_{L^{ p}(B_L(0))}
\| \tilde\phi_j\|_{W^{1,p'}(B_{2L}(0))}.
\end{split}
\end{equation}

Combining  \eqref{n39-29-5}, \eqref{nn39-29-5} and \eqref{30-17-8}, we are led to
\begin{equation}\label{39-29-5}
\begin{split}
&s_{\lambda, j}^2\bar \lambda\int_{B_{2L}(0)}\Bigl(  1_{\{\tilde  {\mathcal U}_{\lambda, {\bf x}_\lambda, {\bf a}_\lambda, j}
 +\tilde \omega_{\lambda, j}>\kappa_{\lambda,j}\}} -1_{\{\tilde U_{\lambda,x_{\lambda,j},a_{\lambda, j}, j}> a_{\lambda, j}\}}
  \Bigr)\tilde\phi_j
  \\
  =&2  \int_{|y|=1}\tilde  \omega_{\lambda, j} \tilde \phi_j
  \\
&
+\Bigl[O\Bigl(\sum_{i=1}^k\frac{s_{\lambda, i}}{|\ln s_{\lambda, i}|}\Bigr)+o(1)\|\omega_{ \lambda}\|_{L^\infty(\Omega)}
+o(1)\|\nabla \tilde \omega_{\lambda, j}\|_{L^{ p}(B_L(0))}
\Bigr]
\|\tilde \phi_j\|_{W^{1,p'}(B_{2L}(0))}.
\end{split}
\end{equation}

 From  \eqref{37-29-5}  and \eqref{39-29-5}, we obtain
\begin{equation}\label{38-29-5}
\begin{split}
&\bigl\langle  R_\lambda(\omega_\lambda), \phi\bigr\rangle\\
=& \Bigl[O\Bigl(\sum_{i=1}^k\frac{s_{\lambda, i}}{|\ln s_{\lambda, i}|}\Bigr)+o(1)\|\omega_{ \lambda}\|_{L^\infty(\Omega)}
+o(1)\|\nabla \tilde \omega_{\lambda, j}\|_{L^{ p}(B_L(0))}
\Bigr]
\|\tilde \phi_j\|_{W^{1,p'}(B_{2L}(0))},
\end{split}
\end{equation}
which implies
\begin{equation}\label{40-29-5}
\begin{split}
&  s_{\lambda, 1}^{\frac2{p'}-1} \| R_\lambda(\omega_\lambda))\|_{W^{-1,p}(\cup_{j=1}^k B_{2Ls_{\lambda,j}}(x_{\lambda, j}))}\\
=&
O\Bigl(\sum_{i=1}^k\frac{s_{\lambda, i}}{|\ln s_{\lambda, i}|}\Bigr)+o(1)\|\omega_{ \lambda}\|_{L^\infty(\Omega)}+o(1)s_{\lambda, j}^{1-\frac2p}\|\nabla
\omega_{\lambda}\|_{L^{p}(\cup_{j=1}^k B_{2Ls_{\lambda,j}}(x_{\lambda, j}))}.
\end{split}
\end{equation}
  So from \eqref{35-29-5}, we derive
\begin{equation}\label{31-17-8}
\begin{split}
&s_{\lambda, 1}^{1-\frac2p}\|\nabla
\omega_\lambda\|_{L^{p}(\cup_{j=1}^k B_{2Ls_{\lambda,j}}(x_{\lambda, j}))}  + \|
\omega_\lambda\|_{L^\infty(\Omega)}
\\
\le  &C\sum_{i=1}^k\frac{s_{\lambda, i}}{|\ln s_{\lambda, i}|}+o(1)\|\omega_{ \lambda}\|_{L^\infty(\Omega)}
+o(1)s_{\lambda, j}^{1-\frac2p}\|\nabla
\omega_{\lambda}\|_{L^{p}(\cup_{j=1}^k B_{2Ls_{\lambda,j}}(x_{\lambda, j}))},
  \end{split}
\end{equation}
which gives the result.

\end{proof}

\section{The Necessary condition for the location of the vortices }

 Using Proposition~\ref{p34}, we can  improve the  estimate for $\Gamma_{\lambda, j}$ in Lemma~\ref{al1} as follows.

 \begin{proposition}\label{p1-17-8}

 The set
\[
\Gamma_{\lambda, j}= \bigl\{ x:\; u_\lambda(s_{\lambda, j}x+ x_{\lambda, j})= \kappa_{\lambda,j}\bigr\}\cap B_{L}(0)
\]
is a continuous closed curve in $\mathbb R^2$, and
\begin{equation}\label{70-17-8}
\Gamma_{\lambda, j}(\theta) = (1+\tilde t_{\lambda, j}(\theta) ) ( \cos \theta, \sin\theta)\\
=( \cos \theta, \sin\theta)+O(s_{\lambda, j}),
\quad \theta\in [0, 2\pi],
\end{equation}
for some function $\tilde t_{\lambda, j}(\theta)$.

\end{proposition}

In the following, we will use $D_x$ to denote the partial derivative for $G(y, x)$ with respect to $x$.

We can also improve Lemma~\ref{l1-28-7} as follows.

\begin{proposition}\label{p2-17-8}

 It holds
\begin{equation}\label{add1}
u_\lambda(x) =\bar \lambda\sum_{j=1}^k |\Omega_{\lambda, j}| G( x_{\lambda, j}, x)
 +  O\Bigl(\sum_{j=1}^k \frac{\bar \lambda s^4_{\lambda, j} }{|x-x_{\lambda,j}|^2  }   \Bigr),
\end{equation}
and
\begin{equation}\label{00add1}
\frac{\partial u_\lambda(x)}{\partial x_i} =\bar \lambda\sum_{j=1}^k |\Omega_{\lambda, j}| D_{x_i} G( x_{\lambda, j}, x)
 +  O\Bigl(\sum_{j=1}^k \frac{\bar \lambda s^4_{\lambda, j} }{|x-x_{\lambda,j}|^3  }   \Bigr),
\end{equation}
for $x\in  \Omega\setminus \cup_{j=1}^k B_{L s_{\lambda, j}}(x_{\lambda, j})$, where $L>0$ is a fixed large  constant.

\end{proposition}

\begin{proof}

 For any   $x\in  \Omega\setminus \cup_{j=1}^k B_{L s_{\lambda, j}}(x_{\lambda, j})$,
\[
\begin{split}
u_\lambda (x) = &\bar \lambda\sum_{j=1}^k \int_{\Omega_{\lambda, j}} G(y, x)\,dy\\
=&\bar \lambda \sum_{j=1}^k|\Omega_{\lambda, j}| G( x_{\lambda, j}, x)+
 \bar \lambda \sum_{j=1}^k \int_{\Omega_{\lambda, j}} \bigl( G(y, x)- G( x_{\lambda, j}, x)\Bigr)\,dy\\
=&\bar \lambda \sum_{j=1}^k|\Omega_{\lambda, j}| G( x_{\lambda, j}, x)+ \bar \lambda\int_{\Omega_{\lambda, j}}\bigl\langle \nabla G( x_{\lambda, j}, x), y-x_{\lambda, j}\bigr\rangle\,dy
 +  O\Bigl(\sum_{j=1}^k \frac{\bar \lambda s^4_{\lambda, j} }{|x-x_{\lambda, j}|^2  }   \Bigr).
\end{split}
\]

By Proposition~\ref{p1-17-8}. we find
\[
\begin{split}
&\int_{\Omega_{\lambda, j}}\bigl\langle \nabla G( x_{\lambda, j}, x), y-x_{\lambda, j}\bigr\rangle\,dy\\
=&\int_{\Omega_{\lambda, j}\setminus  B_{s_{\lambda, j}}(x_{\lambda, j})}\bigl\langle \nabla G( x_{\lambda, j}, x), y-x_{\lambda, j}\bigr\rangle\,dy
+\int_{ B_{s_{\lambda, j}}(x_{\lambda, j})}\bigl\langle \nabla G( x_{\lambda, j}, x), y-x_{\lambda, j}\bigr\rangle\,dy
\\
=& O\bigl( \frac{s_{\lambda, j}}{|x-x_{\lambda, j}|}|\Omega_{\lambda, j}\setminus  B_{s_{\lambda, j}}(x_{\lambda, j})
|\bigr)=O\bigl(\frac{s^4_{\lambda, j}}{|x-x_{\lambda, j}|} \bigr).
\end{split}
\]
So, we prove \eqref{add1}.
 Similarly, we can prove \eqref{00add1}.

\end{proof}

 Let  $l_\lambda= \min_{i\ne j}| x_{\lambda, i}- x_{\lambda, j}|$, $d_\lambda= \min_{i} d(x_{\lambda, i}, \partial\Omega)$  and  $\tau_\lambda= \min
 \{l_\lambda, d_\lambda\}$.  We know that $\frac{\tau_\lambda}{s_{\lambda, j}}\to +\infty$  as $\lambda \to +\infty$.

\begin{proof}[Proof of Theorem~\ref{th11}]

    We have the following Pohozaev identity:
\begin{equation}\label{1-18-8}
 -\int_{\partial B_{\frac14 \tau_{\lambda}}(x_{\lambda, j})}\frac{\partial u_{\lambda}}{\partial \nu}\frac{\partial u_{\lambda}}{\partial x_i}
+\frac12 \int_{\partial B_{\frac14 \tau_{\lambda}}(x_{\lambda, j })}|\nabla
u_{\lambda}|^2 \nu_i =0,
\end{equation}
where  $\nu=(\nu_1,\nu_2)$ is the outward unit normal of $\partial
B_{\frac14 \tau_{\lambda}}(x_{\lambda,j })$.

Using Proposition~\ref{p2-17-8}, we obtain
\begin{equation}\label{2-18-8}
\begin{split}
& -\sum_{l=1}^k\sum_{h=1}^k\int_{\partial B_{\frac14 \tau_{\lambda}}(x_{\lambda,
j})}|\Omega_{\lambda, l}|
|\Omega_{\lambda, h}|\bigl\langle D G( x_{\lambda, h}, x), \nu\bigr\rangle D_{x_i} G( x_{\lambda, l}, x)\\
& +\frac12\int_{\partial B_{\frac14 \tau_{\lambda}}(x_{\lambda, j })}
\Bigl(|\sum_{l=1}^k |\Omega_{\lambda, l}|D G( x_{\lambda,
l}, x)|\Bigr)^2 \nu_i =O\Bigl(  \frac{  s^6_{\lambda, j}}{\tau_\lambda^3}  \Bigr).
\end{split}
\end{equation}

 For any $\tau>\theta>0$ small, it holds
\begin{equation}\label{3-18-8}
 \begin{split}
& -\int_{\partial B_\tau(x_{\lambda, j})}\bigl\langle D G( x_{\lambda, h}, x), \nu\bigr\rangle D_{x_i} G( x_{\lambda, l}, x)
+\frac12 \int_{\partial B_\tau(x_{\lambda.j })}\bigl\langle D G( x_{\lambda, h}, x), D G( x_{\lambda, l}, x)\bigr\rangle \nu_i\\
=&-\int_{\partial B_\theta(x_{\lambda, j})}
\bigl\langle D G( x_{\lambda, h}, x), \nu\bigr\rangle D_{x_i} G( x_{\lambda, l}, x)
+\frac12 \int_{\partial B_\theta(x_{\lambda, j })}\bigl\langle D G(x, x_{\lambda, h}), D G(x_{\lambda, l}, x)\bigr\rangle \nu_i.
\end{split}
\end{equation}

On the other hand, from Proposition~\ref{p1-17-8},  we find
\begin{equation}\label{4-18-8}
 |\Omega_{\lambda, l}|= \pi  s_{\lambda, l}^2+ O(s_{\lambda, l}^3),
\end{equation}
and
\begin{equation}\label{5-18-8}
 s_{\lambda, j}=\frac{ \sqrt{2 \kappa_{j}}}{\lambda}\Bigl(1+O\bigl(\frac{\ln \ln\lambda }{\ln\lambda}\bigr)\Bigr)
 .
\end{equation}
So, we find that \eqref{1-18-8} and \eqref{2-18-8} imply that
\begin{equation}\label{100-1-12-9}
-\sum_{i\neq j}^k(\kappa_i\kappa_j+o(1)) D G(x_{\lambda, i},
x)\bigr|_{x=x_{\lambda, j}}+(\kappa_i^2+o(1)) \nabla \varphi(x_{\lambda, j})=
O\Bigl(  \frac{  s^2_{\lambda, j}}{\tau_\lambda^3}  \Bigr)=o\bigl(\frac{  1}{\tau_\lambda}  \bigr).
\end{equation}

 We claim  that  $\tau_\lambda\ge c_0>0$.
  So, from \eqref{100-1-12-9}£¬ we find  ${\bf x}_0$ is a critical point of
$\mathcal{W}$.

To show our claim, we argue by contradiction.  Suppose that  $\tau_\lambda\to 0$.  We have two cases:  (i) $l_\lambda =o(1) d_\lambda$;
(ii)  $l_\lambda \ge c d_\lambda$ for some $c>0$.

Case (i).  In this case, we have  $\tau_\lambda= l_\lambda$.  From

\[
 D H(x_{\lambda, i},
x)\bigr|_{x=x_{\lambda, j}}= O\bigl(\frac1{d_\lambda}\bigr)=o\bigl(\frac{  1}{\tau_\lambda}  \bigr),
\]
we find from \eqref{100-1-12-9}£¬
\begin{equation}\label{10-21-10}
\sum_{i\neq j}^k(\kappa_i\kappa_j+o(1)) \frac{ x_{\lambda, i}-
x_{\lambda, j}}{| x_{\lambda, i}-
x_{\lambda, j}|^2    }=o\bigl(\frac{  1}{\tau_\lambda}  \bigr), \quad j=1, \cdots, k.
\end{equation}

 There exists a subset $J$ of $\{1, \cdots, k\}$, such that  for any $j_1,\, j_2\in J$, $ j_1\ne j_2$, it holds
 $ |x_{\lambda, j_1}-x_{\lambda, j_2}|\le C l_\lambda$,  and $ \frac{|x_{\lambda, i}-x_{\lambda, j}|}{ l_\lambda}\to +\infty$, $i\notin J$ and $j\in J$.
Then, \eqref{10-21-10} becomes
\begin{equation}\label{11-21-10'}
\sum_{i\neq j, i\in J}(\kappa_i\kappa_j+o(1)) \frac{ x_{\lambda, i}-
x_{\lambda, j}}{| x_{\lambda, i}-
x_{\lambda, j}|^2    }=o\bigl(\frac{  1}{\tau_\lambda}  \bigr), \quad j\in J.
\end{equation}
 We may assume that $| (x_{\lambda, i}-
x_{\lambda, j})_1|\ge  c' l_\lambda$  for some $c'>0$, where $y_1$ is the first coordinate of $y$.  Then, from \eqref{11-21-10'},
\begin{equation}\label{11-21-10}
\sum_{i\neq j, i\in J_1}(\kappa_i\kappa_j+o(1)) \frac{ (x_{\lambda, i}-
x_{\lambda, j})_1}{| x_{\lambda, i}-
x_{\lambda, j}|^2    }=o\bigl(\frac{  1}{\tau_\lambda}  \bigr), \quad j\in J_1
\end{equation}
for some subset of $J_1$ of $J$, satisfying $|(x_{\lambda, i}-
x_{\lambda, j})_1|\ge c'' l_\lambda $ for any $i,\, j\in J_1$ and $i\ne j$.
Clearly, \eqref{11-21-10} is not true at $x_{\lambda,j}$ with $(x_{\lambda,j})_1= \max_{i\in J_1} (x_{\lambda,i})_1$.

Case (ii). In this case, $d_\lambda\to 0$.  Take $x_{\lambda, j}$ satisfying $d(x_{\lambda, j}, \partial \Omega) = d_\lambda$.
Now we consider \eqref{100-1-12-9}  at $x_{\lambda, j}$. It holds
\[
\frac{\partial  \varphi(x_{\lambda, j})}{\partial \nu}=\frac{\alpha}{d_\lambda} (1+o(1)),
\]
where $\nu$ is the outward unit normal of $\partial\Omega$ at $\bar x_{\lambda, j}$ with $|x_{\lambda, j}-\bar x_{\lambda, j}|=d_\lambda$,
and $\alpha>0$ is some constant.

On the other hand,   if $|x_{\lambda, j}- x_{\lambda, i}|\ge L  d_\lambda$, then
\[
 D G(x_{\lambda, i}, x)\bigr|_{x=x_{\lambda, j}}
= O\bigl(\frac1{L d_\lambda}\bigr).
\]

If $|x_{\lambda, j}- x_{\lambda, i}|\le L  d_\lambda$, then it is easy to check
\[
-\bigl\langle  D G(x_{\lambda, i}, x)\bigr|_{x=x_{\lambda, j}}, \nu\bigr\rangle
\ge  \frac{c''}{|x_{\lambda, j}- x_{\lambda, i}|}\bigl\langle \frac{x_{\lambda, j}- x_{\lambda, i}}{|x_{\lambda, j}- x_{\lambda, i}|},
\nu\bigr\rangle\ge  \frac{o(1)}{d_\lambda},
\]
since $|x_{\lambda, j}- x_{\lambda, i}|\ge c'  d_\lambda$.
So, from \eqref{100-1-12-9},
we find
\[
\frac1{d_\lambda}(1-o(1)) \le o\bigl(\frac1{d_\lambda}\bigr).
\]
This is a contradiction, which concludes the proof of Proposition~\ref{p2-17-8}.

\end{proof}

\begin{remark}\label{re1-18-8}

 If $k=1$, from
 \eqref{2-18-8}, we find
 \begin{equation}\label{6-18-8}
 -\int_{\partial B_{\frac14 \tau_\lambda}(x_{\lambda})}\bigl\langle D G( x_{\lambda, h}, x), \nu\bigr\rangle D_{x_i} G( x_{\lambda, l}, x)
+\frac12\int_{\partial B_{\frac14 \tau_\lambda}(x_{\lambda})}
|D G( x_{\lambda}, x)|^2 \nu_i
=O(s_{\lambda}^2),
\end{equation}
which gives
\begin{equation}\label{16-18-8}
 \nabla H(x_\lambda,  x_\lambda)=O(s_{\lambda}^2).
\end{equation}

\end{remark}

\section{Uniqueness  Results}

In this section, we will prove that if $\Omega$ is convex, the vortex patch problem has unique solution.

Suppose that $\Omega$ is convex.  Then, if $k\ge 2$,  $\mathcal{W}$ has no critical point. Thus, \eqref{1'} has no solution
for $k\ge 2$.
 So, we only need to
consider the case $k=1$ in \eqref{1'}.
Note that if $k=1$, from \cite{CF1},  $H(x,x)$ has a unique critical point $x_0$, which is also non-degenerate.

In this section, we will  prove the following local uniqueness result stated in Theorem~\ref{th1-21-10}.
We argue by contradiction.
Suppose that \eqref{1'}-\eqref{2'}  has two different solutions $u^{(1)}_{\lambda}  $  and
$u^{(2)}_{\lambda} $, which blow up at $x_0$.  We will use
$x_{\lambda}^{(i)}$, $s_{\lambda}^{(i)}$ and $\kappa_{\lambda}^{(i)}$ to denote the parameters appearing in
 $u^{(i)}_{^{\lambda}}  $.

 Let
\begin{equation}\label{1-7-4}
\xi_\lambda(y) =\frac{u^{(1)}_{\lambda}(y)-u^{(2)}_{\lambda}(y)}{\|u^{(1)}_{\lambda}-u^{(2)}_{\lambda}\|_{L^\infty(\Omega)}}.
\end{equation}
Then, $\xi_\lambda$ satisfies $\|\xi_\lambda\|_{L^\infty(\Omega)}=1$ and
\begin{equation}\label{2-7-4}
-\Delta \xi_\lambda = f_\lambda(y),
\end{equation}
where
\begin{equation}\label{3-7-4}
f_\lambda(y)=\frac{\bar \lambda}
{\|u^{(1)}_{\lambda}-u^{(2)}_{\lambda}\|_{L^\infty(\Omega)}}
 \bigl( 1_{B_\delta(x_{0})}
1_{ \{ u^{(1)}_{\lambda} >\kappa_\lambda^{(1)}\}}- 1_{B_\delta(x_{0})}
1_{ \{ u^{(2)}_{\lambda} >\kappa_\lambda^{(2)}\}}\bigr).
\end{equation}

Using the non-degeneracy of $x_0$ and \eqref{16-18-8}, we find $|x^{(i)}_\lambda -x_0|= O(s_\lambda^2)$, $i=1, 2$. This results in
 $f_\lambda(y)=0$ in $\Omega\setminus B_{L s^{(1)}_{\lambda}}(x_\lambda^{(1)})$. So it holds  $\Delta \xi_\lambda=0$ in $\Omega\setminus B_{L s^{(1)}_{\lambda}}(x_\lambda^{(1)})$. To obtain
a contradiction, we only need to prove $\xi_\lambda=o(1)$ in $ B_{L s^{(1)}_{\lambda}}(x_\lambda^{(1)})$.

 For simplicity, in the sequel we use $s_\lambda$ to denote $s^{(1)}_{\lambda}$.   Let
$\Omega_\lambda=\{ y\in \mathbb R^2:\, s_\lambda y+x_\lambda^{(1)}\in \Omega\}.$

\begin{lemma}\label{l10-18-8}

It holds  $\|s_{\lambda}^2f_\lambda(s_{\lambda} y+ x_\lambda^{(1)})\|_{W^{-1,p}  (\Omega_\lambda) }\le C$.  As a result,
\[
 \xi_\lambda(s_{\lambda} y+ x_\lambda^{(1)}) \to \xi,\quad \text{in}\quad  C_{loc}(\mathbb R^2),
\]
as $\lambda\to +\infty$.  Moreover,
\[
\int_{ \Omega_\lambda } s_{\lambda}^2f_\lambda(s_{\lambda} y+ x_\lambda^{(1)})\phi \to  2\int_{|y|=1}\xi\phi,\quad  \forall \,\,\phi\in C_0^\infty(\mathbb R^2).
\]

\end{lemma}

\begin{proof}
  Let
\[
  \tilde \Gamma_\lambda ^{(i)}=\{ y:  \;  u_\lambda^{(i)}( s_\lambda y+ x_{\lambda}^{(1)})= \kappa_\lambda^{(i)}\}.
  \]

 For any $y_\lambda^{(2)}\in\tilde  \Gamma_\lambda ^{(2)}$, let  $ y_\lambda^{(1)}= (1+ t_\lambda) y_\lambda^{(2)}\in \tilde \Gamma_\lambda ^{(1)}$. Then $t_\lambda \to 0$.
 We have
\[
 \begin{split}
 &u_\lambda^{(1)}(s_{\lambda} y_\lambda^{(2)}+ x_{\lambda}^{(1)})-u_\lambda^{(2)}(s_{\lambda} y_\lambda^{(2)}+ x_{\lambda}^{(1)})=u_\lambda^{(1)}(s_{\lambda} y_\lambda^{(2)}+ x_{\lambda}^{(1)})-u_\lambda^{(1)}(s_{\lambda} y_\lambda^{(1)}+ x_{\lambda}^{(1)})\\
 =&- \frac{\partial }{\partial r}u_\lambda^{(1)}(s_{\lambda}(1+ \theta t_\lambda) y_\lambda^{(2)}+ x_{\lambda}^{(1)}) t_\lambda |y_\lambda^{(2)}|
 =\bar \lambda s_{\lambda}^2\bigl(\frac12+ o(1)\bigr) t_\lambda |y_\lambda^{(2)}|,
 \end{split}
 \]
 which, together with  $|y_\lambda^{(2)}|=1+o(1)$, implies
 \begin{equation}\label{20-18-8}
 t_\lambda = \frac{1}{\bar \lambda s_{\lambda}^2}\bigl(2+o(1)\bigr)\Bigl(u_\lambda^{(1)}(s_{\lambda} y_\lambda^{(2)}+ x_{\lambda}^{(1)})-u_\lambda^{(2)}(s_{\lambda} y_\lambda^{(2)}+ x_{\lambda}^{(1)})\Bigr).
 \end{equation}
 As a result,
\begin{equation}\label{21-18-8}
 \begin{split}
&\int_{ \Omega_\lambda}  s_{\lambda}^2f_\lambda(s_{\lambda} y+ x_\lambda^{(1)})\phi \\
=& \frac{ \bar \lambda s_{\lambda}^2 }{\|u^{(1)}_{\lambda}-u^{(2)}_{\lambda}\|_{L^\infty(\Omega)}}\int_{y\in\tilde \Gamma^{(1)}_\lambda}\int_0^{t_\lambda}\phi((1+t)y)\,dt\,dy \\
=& \frac{ \bar \lambda s_{\lambda}^2 }{\|u^{(1)}_{\lambda}-u^{(2)}_{\lambda}\|_{L^\infty(\Omega)}}\int_{y\in\tilde \Gamma^{(1)}_\lambda}
\phi(y)t_\lambda\,dy \\
&+\frac{ \bar \lambda s_{\lambda}^2 }{\|u^{(1)}_{\lambda}-u^{(2)}_{\lambda}\|_{L^\infty(\Omega)}}\int_{y\in\tilde \Gamma^{(1)}_\lambda}
\int_0^{t_\lambda}[\phi(y+ty)-\phi(y)]\,dy\,dt \\
=& O\Bigl( \int_{y\in\tilde \Gamma^{(1)}_\lambda}|\phi(y)|\,dy\Bigl)+\frac{ \bar \lambda s_{\lambda}^2 }{\|u^{(1)}_{\lambda}-u^{(2)}_{\lambda}\|_{L^\infty(\Omega)}}\int_{y\in\tilde \Gamma^{(1)}_\lambda}
\int_0^{t_\lambda} \bigl\langle \nabla\phi(y+\theta ty), yt\bigr\rangle \,dt\,dy\\
= & O\Bigl( \int_{y\in\tilde \Gamma^{(1)}_\lambda}|\phi(y)|\,dy
+\int_0^{t_\lambda} \int_{y\in\tilde \Gamma^{(1)}_\lambda}|\nabla\phi(y+\theta ty)|\,dy dt\Bigr)\\
=& O\Bigl( \|\phi\|_{W^{1,p'}(\mathbb R^2)}\Bigl).
\end{split}
 \end{equation}
 So, $\| s_{\lambda}^2f_\lambda(s_{\lambda} y+ x_\lambda^{(1)})
 \|_{W^{-1,p}(\Omega_\lambda)}\le C$.

 Similar to \eqref{21-18-8}, using \eqref{20-18-8}, we can prove
\begin{equation}\label{22-18-8}
 \begin{split}
&\int_{ \Omega_\lambda}  s_{\lambda}^2f_\lambda(s_{\lambda} y+ x_\lambda^{(1)})\phi
 \\
=& \frac{ \bar \lambda s_{\lambda}^2 }{\|u^{(1)}_{\lambda}-u^{(2)}_{\lambda}\|_{L^\infty(\Omega)}}\int_{y\in\tilde \Gamma^{(1)}_\lambda}
\phi(y)t_\lambda\,dy
+o(1)\\
=&\int_{y\in\tilde \Gamma^{(1)}_\lambda}(2+o(1))\xi_\lambda (s_{\lambda} y^{(2)}+ x_{\lambda}^{(1)}) \phi(y)\to  2\int_{|y|=1}\xi\phi.
 \end{split}
 \end{equation}

\end{proof}

\begin{lemma}\label{l1-10-4}

 It holds
\begin{equation}\label{10-9-4}
\tilde \xi_{\lambda}(y):= \xi_{\lambda}(s_\lambda y+x_\lambda^{(1)}) \to b_{1}\frac{\partial w}{\partial x_1}+ b_{2}\frac{\partial w}{\partial x_2},
\end{equation}
uniformly in $C(B_R(0))$ for any $R>0$, where $b_{1}$ and $b_2$
are some constants.

\end{lemma}

\begin{proof}

It follows from Lemma~\ref{l10-18-8} that $\xi$ satisfies
\begin{equation}\label{3-8-4}
-\Delta \xi = 2\delta_{|y|=1}\xi, \quad \text{in}\;  \mathbb R^2.
\end{equation}
So the result follows.

\end{proof}

 To prove $\xi_\lambda=o(1)$ in $ B_{Ls_{\lambda}}(x_\lambda^{(1)})$, we need to prove $b_1= b_2=0$.  We will use a local
 Pohozaev identity to achieve this.  The following lemma gives an estimate for $x$  in $\Omega \setminus   B_{2\theta}(x^{(1)}_{\lambda})$.

\begin{lemma}\label{l1-9-4}

We have the following estimate:
\begin{equation}\label{1-9-4}
\begin{split}
  \xi_\lambda(x) =&   A_{\lambda} G(x^{(1)}_{\lambda}, x)
  + \sum_{h=1}^2 B_{ h, \lambda} \partial _{h}G(x^{(1)}_{\lambda}, x)\\
&
+O\bigl( s_\lambda^2\bigr),\quad
\text{in $C^1\bigl(\Omega \setminus   B_{2\theta}(x^{(1)}_{\lambda})\bigr)$},
\end{split}
\end{equation}
where $\theta>0$ is any   small constant, $\partial_h G(y, x)=\frac{\partial G(y,x)}{\partial y_h}$,
\begin{equation}\label{10-15-5-14}
A_{\lambda}= \int_{\Omega}  f_\lambda(y) \,dy,
\end{equation}
and
\begin{equation}\label{nn1-9-4}
 B_{ h, \lambda}= \int_{ \Omega } ( y_h- x^{(1)}_{\lambda, h})
f_\lambda(y )\,dy.
\end{equation}

\end{lemma}

\begin{proof}

   We have
\begin{equation}\label{02-10-4}
\begin{split}
&\xi_\lambda(x)
= \int_{\Omega} G(y,x) f_\lambda (y)\,dy\\
=&  A_{\lambda} G(x^{(1)}_{\lambda}, x)+  \sum_{h=1}^2 B_{ h, \lambda} \partial _{h}G(x^{(1)}_{\lambda}, x)\\
 &+
  \int_{\Omega}\bigl( G(y,x) -G(x^{(1)}_{\lambda}, x)- \langle \nabla G(x^{(1)}_{\lambda}, x), y-x_\lambda^{(1)}\rangle\bigr) f_\lambda(y)\,dy.
 \end{split}
\end{equation}

 On the other hand, similar to \eqref{21-18-8}, we can prove
 \[
\int_{\Omega}| f_\lambda(y)|\,dy= \int_{ \Omega_\lambda}  s_{\lambda}^2| f_\lambda(s_{\lambda} y+ x_\lambda^{(1)})|
=O(1).
\]
So
\begin{equation}\label{00-02-10-4}
\begin{split}
& \bigl|\int_{\Omega}\bigl( G(y,x) -G(x^{(1)}_{\lambda}, x)- \langle \nabla G(x^{(1)}_{\lambda}, x), y-x_\lambda^{(1)}\rangle\bigr) f_\lambda(y)\,dy
\bigr|\\
\le & C s_\lambda^2 \int_{\Omega}| f_\lambda(y)|\,dy\le C s_\lambda^2.
\end{split}
\end{equation}

Similarly, we can prove that \eqref{1-9-4} holds in $C^1\bigl(\Omega \setminus   B_{2\theta}(x^{(1)}_{\lambda})\bigr)$.

\end{proof}

In the following, we will use $\partial $ or $\nabla $ to denote the partial derivative for any function $h(y, x)$ with respect to $y$, while
we will use $D $ to denote the partial derivative for any function $h(y, x)$ with respect to $x$.

\begin{proof}[Proof of Theorem~\ref{th1-21-10}]

We have the following Pohozaev identity for $\xi_\lambda$:  For $0<d<\delta$, it holds
\begin{equation}\label{14-11-9}
 \begin{split}
&-\int_{\partial B_d(x^{(1)}_{\lambda})}\frac{\partial  \xi_{\lambda}}{\partial \nu}\frac{\partial  u^{(1)}_{\lambda}}
{\partial x_i}-\int_{\partial B_d(x^{(1)}_{\lambda})}\frac{\partial  u^{(2)}_{\lambda}}{\partial \nu}\frac{\partial  \xi_{\lambda}}{\partial x_i}\\
&+\frac12 \int_{\partial B_d(x^{(1)}_{\lambda})}\bigl\langle \nabla ( u^{(1)}_{\lambda} + u^{(2)}_{\lambda}),
\nabla  \xi_{\lambda}\bigr\rangle \nu_i=0.
\end{split}
\end{equation}

By Proposition~\ref{p2-17-8} and Lemma~\ref{l1-9-4}, we obtain from \eqref{14-11-9}
\begin{equation}\label{1-19-8}
 \begin{split}
&-\int_{\partial B_d(x^{(1)}_{\lambda})}\frac{\partial  \xi_{\lambda}}{\partial \nu} D_{x_i}  G( x_\lambda^{(1)}, x)-\int_{\partial B_d(x^{(1)}_{\lambda})}\bigl\langle
  G( x^{(1)}_{\lambda}, x),  \nu\bigr\rangle \frac{\partial  \xi_{\lambda}}{\partial x_i}\\
&+ \int_{\partial B_d(x^{(1)}_{\lambda})}\bigl\langle D  G( x^{(1)}_{\lambda}, x),
\nabla  \xi_{\lambda}\bigr\rangle \nu_i=O\bigl(s_\lambda^2\bigr).
\end{split}
\end{equation}

On the other hand, from Lemma~\ref{l1-9-4}  and Remark~\ref{re1-18-8},
 \eqref{1-19-8}  becomes
\begin{equation}\label{3-19-8}
 \begin{split}
&-\int_{\partial B_d(x^{(1)}_{\lambda})} \sum_{h=1}^2 B_{ h, \lambda}\bigl\langle D  \partial _{h}G(x^{(1)}_{\lambda}, x),  \nu\bigr\rangle D_{x_i}  G( x_\lambda^{(1)}, x)
\\
&-\int_{\partial B_d(x^{(1)}_{\lambda})}\sum_{h=1}^2 B_{ h, \lambda}\bigl\langle
  G( x^{(1)}_{\lambda}, x),  \nu\bigr\rangle
  D_{x_i}  \partial _{h}G(x^{(1)}_{\lambda}, x)\\
&+ \int_{\partial B_d(x^{(1)}_{\lambda})} \sum_{h=1}^2 B_{ h, \lambda}\bigl\langle D  G( x^{(1)}_{\lambda}, x),
D \partial _{h}G(x^{(1)}_{\lambda}, x)\bigr\rangle \nu_i=O\bigl(s_\lambda^2\bigr).
\end{split}
\end{equation}

We can check for any small $\theta>0$,
\begin{equation}\label{4-19-8}
 \begin{split}
&-\int_{\partial B_d(x^{(1)}_{\lambda})} \bigl\langle D  \partial _{h}G(x^{(1)}_{\lambda}, x),  \nu\bigr\rangle D_{x_i}  G( x_\lambda^{(1)}, x)
-\int_{\partial B_d(x^{(1)}_{\lambda})}\bigl\langle
  G( x^{(1)}_{\lambda}, x),  \nu\bigr\rangle
  D_{x_i}  \partial _{h}G(x^{(1)}_{\lambda}, x)\\
&+ \int_{\partial B_d(x^{(1)}_{\lambda})}\bigl\langle D  G( x^{(1)}_{\lambda}, x),
D \partial _{h}G(x^{(1)}_{\lambda}, x)\bigr\rangle \nu_i
\\
=&-\int_{\partial B_\theta(x^{(1)}_{\lambda})} \bigl\langle D  \partial _{h}G(x^{(1)}_{\lambda}, x),  \nu\bigr\rangle D_{x_i}  G( x_\lambda^{(1)}, x)
-\int_{\partial B_\theta(x^{(1)}_{\lambda})}\bigl\langle
  G( x^{(1)}_{\lambda}, x),  \nu\bigr\rangle
  D_{x_i}  \partial _{h}G(x^{(1)}_{\lambda}, x)\\
&+ \int_{\partial B_\theta(x^{(1)}_{\lambda})}\bigl\langle D  G( x^{(1)}_{\lambda}, x),
D \partial _{h}G(x^{(1)}_{\lambda}, x)\bigr\rangle \nu_i.
\end{split}
\end{equation}

We define the following quadric form
\begin{equation}\label{5-19-8}
Q(u, v)
=-\int_{\partial B_\theta(x^{(1)}_{\lambda})} \frac{\partial  v}{\partial \nu} \frac{\partial u}{\partial x_i}-\int_{\partial B_\theta(x^{(1)}_{\lambda})}\frac{\partial u
}{\partial \nu}\frac{\partial  v}{\partial x_i}
+ \int_{\partial B_\theta(x^{(1)}_{\lambda})} \bigl\langle \nabla u,
\nabla v\bigr\rangle \nu_i.
\end{equation}
Note that if $u$ and $v$ are harmonic in $ B_\delta(x^{(1)}_{\lambda})\setminus \{x^{(1)}_{\lambda}\}$, then $Q(u, v)$ is independent of $\theta>0$.

Denote $G(y, x) = S(y, x) -H(y,x )$ and $S(y, x)=\frac1{2\pi}\ln\frac1{|y-x|}$ is the singular part of $G(y, x)$.
Note that
$Q(S(x^{(1)}_{\lambda}, x),  \partial _{h}S(x^{(1)}_{\lambda}, x))$ is either infinity or zero since the singularity involved is
of order $\frac1{|x-x_\lambda|^3}$. On the other hand, $Q(G(x^{(1)}_{\lambda}, x),  \partial _{h}G(x^{(1)}_{\lambda}, x))$ is bounded if  $\theta> 0$ is fixed.
Thus $Q(S(x^{(1)}_{\lambda}, x),  \partial _{h}S(x^{(1)}_{\lambda}, x))=0$.
Therefore,
\begin{equation}\label{6-19-8}
 \begin{split}
&Q(G(x^{(1)}_{\lambda}, x),  \partial _{h}G(x^{(1)}_{\lambda}, x))\\
=&Q(S(x^{(1)}_{\lambda}, x),  \partial _{h}H(x^{(1)}_{\lambda}, x))+Q(H(x^{(1)}_{\lambda}, x),  \partial _{h}S(x^{(1)}_{\lambda}, x))+o_\theta(1).
\end{split}
\end{equation}

Direct calculations show that
\begin{equation}\label{10-19-8}
 Q(S(x^{(1)}_{\lambda}, x),  \partial _{h}H(x^{(1)}_{\lambda}, x))\\
= D_{x_i} \partial_h H(x^{(1)}_{\lambda}, x^{(1)}_{\lambda})+o_\theta(1),
\end{equation}
and
\begin{equation}\label{12-19-8}
 \begin{split}
&Q(H(x^{(1)}_{\lambda}, x),  \partial _{h}S(x^{(1)}_{\lambda}, x))\\
=&-\int_{\partial B_\theta(x^{(1)}_{\lambda})} \bigl\langle D  \partial _{h}S(x^{(1)}_{\lambda}, x),  \nu\bigr\rangle \Bigl( D_{x_i}  H( x^{(1)}_{\lambda}, x)
-D_{x_i}  H( x^{(1)}_{\lambda}, x^{(1)}_{\lambda})\Bigr)\\
&-\int_{\partial B_\theta(x^{(1)}_{\lambda})}\Bigl(\bigl\langle D H( x^{(1)}_{\lambda}, x), \nu\bigr\rangle-
\bigl\langle D  H(x^{(1)}_{\lambda}, x^{(1)}_{\lambda}), \nu\bigr\rangle\Bigr) D_{x_i}  \partial _{h}S(x^{(1)}_{\lambda}, x)
\\
&+ \int_{\partial B_\theta(x^{(1)}_{\lambda})} \bigl\langle D  H( x^{(1)}_{\lambda}, x)- D H(x^{(1)}_{\lambda}, x^{(1)}_{\lambda})   ,
D  \partial _{h}S(x^{(1)}_{\lambda}, x)\bigr\rangle \nu_i\\
=&-\int_{\partial B_\theta(x^{(1)}_{\lambda})}\bigl\langle D  \partial _{h}S(x^{(1)}_{\lambda}, x),  \nu\bigr\rangle
  \bigl\langle D D_{x_i}  H(x^{(1)}_{\lambda}, x^{(1)}_{\lambda}), x-x^{(1)}_{\lambda}\bigr\rangle\\
&-\int_{\partial B_\theta(x^{(1)}_{\lambda})}\bigl\langle D^2  H(x^{(1)}_{\lambda}, x^{(1)}_{\lambda})(x-x^{(1)}_{\lambda}), \nu\bigr\rangle
\frac{\partial  \partial _{h}S(x^{(1)}_{\lambda}, x)}{\partial x_i}
\\
&+ \int_{\partial B_\theta(x^{(1)}_{\lambda})} \bigl\langle    D^2   H(x^{(1)}_{\lambda}, x^{(1)}_{\lambda})(x-x^{(1)}_\lambda)  ,
D  \partial _{h}S(x^{(1)}_{\lambda}, x)\bigr\rangle \nu_i+o_\theta(1)\\
=&J_1+J_2+ J_3+o_\theta(1).
\end{split}
\end{equation}

  Let   $y= x-x^{(1)}_{\lambda}$.  We have
\[
 \partial _{h}S(x^{(1)}_{\lambda}, x)=\frac1{2\pi} \frac { y_h}{|y|^2},
\]
\[
\bigl\langle D  \partial _{h}S(x^{(1)}_{\lambda}, x),  \nu\bigr\rangle =-\frac1{2\pi} \frac { y_h}{|y|^3},
\]
and
\[
D_{x_t} \partial _{h}S(x^{(1)}_{\lambda}, x)=\frac1{2\pi} \Bigl(\frac { \delta_{ht}}{|y|^2}- \frac {2 y_t y_h}{|y|^4}\Bigr).
\]
So we find
\begin{equation}\label{1-21-8}
J_1= \frac1{2\pi}\int_{|y|=\theta} \frac { y_h}{|y|^3} D^2_{x_h x_i}  H(x^{(1)}_{\lambda}, x^{(1)}_{\lambda}) y_h=\frac12 D^2_{x_h x_i}  H(x^{(1)}_{\lambda}, x^{(1)}_{\lambda}).
\end{equation}

On the other hand,
\begin{equation}\label{2-21-8}
 \begin{split}
&J_2+ J_3\\
=&
-\frac1{2\pi}\int_{ |y|=\theta  }\sum_{l=1}^2\sum_{t=1}^2 D^2_{x_l x_t}  H(x^{(1)}_{\lambda}, x^{(1)}_{\lambda})
\frac{y_l y_t}{|y|}
 \Bigl(\frac { \delta_{hi}}{|y|^2}- \frac {2 y_i y_h}{|y|^4}\Bigr)
\\
&+ \frac1{2\pi}\int_{ |y|=\theta }\sum_{l=1}^2\sum_{t=1}^2 D^2_{x_l x_t}  H(x^{(1)}_{\lambda}, x^{(1)}_{\lambda})
y_l
 \Bigl(\frac { \delta_{ht}}{|y|^2}- \frac {2 y_t y_h}{|y|^4}\Bigr)
 \frac{y_i}{|y|}\\
 =&
-\frac1{2\pi}\int_{ |y|=\theta }\sum_{l=1}^2\sum_{t=1}^2 D^2_{x_l x_t}  H(x^{(1)}_{\lambda}, x^{(1)}_{\lambda})
\frac{\delta_{hi}y_l y_t}{|y|^3}
\\
&+ \frac1{2\pi}\int_{|y|=\theta  }\sum_{l=1}^2\sum_{t=1}^2 D^2_{x_l x_t}  H(x^{(1)}_{\lambda}, x^{(1)}_{\lambda})
 \frac{\delta_{ht} y_ly_i}{|y|^3}.
\end{split}
\end{equation}

If $i\ne h$,  then
\begin{equation}\label{3-21-8}
 \begin{split}
J_2+ J_3&
=
 \frac1{2\pi}\int_{ |y|=\theta }\sum_{l=1}^2\sum_{t=1}^2 D^2_{x_l x_t}  H(x^{(1)}_{\lambda}, x^{(1)}_{\lambda})
 \frac{\delta_{ht} y_ly_i}{|y|^3}\\
 &=\frac1{2\pi}\int_{ |y|=\theta } D^2_{x_i x_h}  H(x^{(1)}_{\lambda}, x^{(1)}_{\lambda})
 \frac{y_i^2}{|y|^3}=\frac12 D^2_{x_i x_h}  H(x^{(1)}_{\lambda}, x^{(1)}_{\lambda}).
\end{split}
\end{equation}

Combining \eqref{1-21-8}  and \eqref{3-21-8}, we obtain
\begin{equation}\label{4-21-8}
J_1+J_2+ J_3
= D^2_{x_i x_h}  H(x^{(1)}_{\lambda}, x^{(1)}_{\lambda}), \quad  i\ne h.
\end{equation}

If $i= h$,  then
\begin{equation}\label{5-21-8}
 \begin{split}
J_2+ J_3
=&
-\frac1{2\pi}\int_{ |y|=\theta }\sum_{l=1}^2\sum_{t=1}^2 D^2_{x_l x_t}  H(x^{(1)}_{\lambda}, x^{(1)}_{\lambda})
\frac{y_l y_t}{|y|^3}
\\
&+ \frac1{2\pi}\int_{ |y|=\theta }\sum_{l=1}^2\sum_{t=1}^2 D^2_{x_l x_t}  H(x^{(1)}_{\lambda}, x^{(1)}_{\lambda})
 \frac{\delta_{it} y_ly_i}{|y|^3}\\
=&  -\frac12  \sum_{l=1}^2 D_{x_l x_l}  H(x^{(1)}_{\lambda}, x^{(1)}_{\lambda})+\frac12 D_{x_i x_i}  H(x^{(1)}_{\lambda}, x^{(1)}_{\lambda})
=\frac12 D_{x_i x_i}  H(x^{(1)}_{\lambda}, x^{(1)}_{\lambda}),
\end{split}
\end{equation}
since
\[
\sum_{l=1}^2 D_{x_l x_l}  H(x^{(1)}_{\lambda}, x)= \Delta_{x} H(x^{(1)}_{\lambda}, x)=0,\quad \forall x\in\Omega.
\]
So, we also have
\begin{equation}\label{6-21-8}
J_1+J_2+ J_3
=
D_{x_i x_i}  H(x^{(1)}_{\lambda}, x^{(1)}_{\lambda}), \quad  i=h.
\end{equation}

Combining \eqref{6-19-8}, \eqref{10-19-8}, \eqref{4-21-8} and \eqref{6-21-8}, we are led to
\begin{equation}\label{7-21-8}
 \begin{split}
&Q(G(x^{(1)}_{\lambda}, x),  \partial _{h}G(x^{(1)}_{\lambda}, x))\\
=&D_{x_i} \partial_h H(x^{(1)}_{\lambda}, x^{(1)}_{\lambda})+D^2_{x_i x_h}  H(x^{(1)}_{\lambda}, x^{(1)}_{\lambda})
=\frac{\partial^2 \varphi(x_\lambda^{(1)})}
{\partial x_i\partial x_h}.
\end{split}
\end{equation}

Thus, \eqref{3-19-8} becomes
\begin{equation}\label{8-21-8}
 (D^2\varphi(x_\lambda^{(1)})) B_\lambda = O(s_{\lambda}^2),
\end{equation}
which, together  with the    non-degeneracy of the  critical point  $x_0$, implies  $B_\lambda =O(s_{\lambda}^2)$.  But
\begin{equation}\label{9-21-8}
 B_{ h, \lambda}= s_\lambda\int_{ \mathbb R^2 } y_h s_\lambda^2f_\lambda(s_\lambda y+x_\lambda^{(1)} )\,dy= s_\lambda \Bigl( \int_{|y|=1}\bigl(   b_{1}\frac{\partial w}{\partial y_1}+ b_{2}\frac{\partial w}{\partial y_2} \bigr) y_h
+o(1)\Bigr).
\end{equation}
Thus, $b_1=b_2=0$. So we have proved  $\xi_\lambda=o(1)$ in $ B_{L s^{(1)}_{\lambda}}(x_\lambda^{(1)})$.  On the other hand,
$\Delta \xi_\lambda=0$ in $ \Omega\setminus B_{L s^{(1)}_{\lambda}}(x_\lambda^{(1)})$ and $\xi_\lambda=0$ on $\partial\Omega$. By
the maximum principle, we conclude $\xi_\lambda=o(1)$ in $\Omega$. This is a contradiction to $\|\xi_\lambda\|_{L^\infty(\Omega)}=1$.

\end{proof}

\appendix

\section{Some essential estimates}

  Let $\mathcal U_{\lambda, {\bf x}_\lambda, {\bf a}}$ be defined in
 \eqref{1-5-11}.  We have

 \begin{lemma}\label{l1-20-10}

 For any fixed $L>0$, it holds
\begin{equation}\label{1-11-9}
\begin{split}
& \mathcal U_{\lambda, {\bf x}_\lambda, {\bf a}_\lambda}(y)-\kappa_{\lambda,i} \\
=&
 U_{\lambda,x_{\lambda,i}, a_{\lambda,i}}(y)- a_{\lambda,i} -
 \frac {a_{\lambda,i}}{\ln \frac{R}{s_{\lambda, i}}}
 \bigl\langle D g(x_{\lambda,i},x_{\lambda,i}), y-x_{\lambda,i}\bigr\rangle
 \\
 &+\sum_{j\ne i}\frac{a_{\lambda,j}}{\ln\frac R{s_{\lambda, j}}
}
\bigl\langle D  \bar G(x_{\lambda,i},x_{\lambda,j}) ,y-x_{\lambda,i}\bigr\rangle
 +O\Bigl( \frac{s_{\lambda,i}^2 }{\ln\lambda } \Bigr),\quad y\in B_{L s_{\lambda,i}}(x_{\lambda,i}).
\end{split}
\end{equation}
\end{lemma}

\begin{proof}

 For $y\in B_{L s_{\lambda,i}}(x_{\lambda,i})$, where $L>0$ is any
fixed constant,
\[
\begin{split}
& PU_{\lambda,x_{ \lambda,i},a_{\lambda, i}}(y)-\kappa_{\lambda,i} = U_{\lambda,x_{\lambda,i}, a_{\lambda,i}}(y)-\kappa_{\lambda,i}  -\frac
{a_{\lambda,i}}{\ln \frac{R}{s_{\lambda, i}}} g(y,x_{\lambda, i})
\\
=&   U_{\lambda,x_{\lambda,i}, a_{\lambda,i}}(y)-\kappa_{\lambda,i}  -\frac
{a_{\lambda,i}}{\ln \frac{R}{s_{\lambda, i}}} g(x_{\lambda, i},x_{\lambda, i})-\frac
{a_{\lambda,i}}{\ln \frac{R}{s_{\lambda, i}}}\Bigl(
 \bigl\langle D g(x_{\lambda,i},x_{\lambda,i}), y-x_{\lambda,i}\bigr\rangle\\
 & +O( |y-x_{\lambda,i}|^2)\Bigr)\\
=&   U_{\lambda,x_{\lambda,i}, a_{\lambda,i}}(y)-\kappa_{\lambda,i} -\frac
{a_{\lambda,i}}{\ln \frac{R}{s_{\lambda, i}}} g(x_{\lambda, i},x_{\lambda, i})-\frac
{a_{\lambda,i}}{\ln \frac{R}{s_{\lambda, i}}}
 \bigl\langle D g(x_{\lambda,i},x_{\lambda,i}), y-x_{\lambda,i}\bigr\rangle+O\bigl(\frac{ s_{\lambda,i}^2}{\ln\lambda}\bigr),
\end{split}
\]
and for $j\ne i$ and  $y\in B_{L s_{\lambda,i}}(x_{\lambda,i})$, by \eqref{1-5-1}
\[
\begin{split}
 PU_{\lambda,x_{\lambda, j},a_{\lambda, j}}(y)=&U_{\lambda, x_{\lambda,j},a_{\lambda,j}}(y)-\frac{a_{\lambda,j}}{\ln\frac R{s_{\lambda, j}}
}
  g(y,x_{\lambda,j})=
 \frac{a_{\lambda,j}}{\ln\frac R{s_{\lambda, j}}
}   \bar G(y,x_{\lambda,j})\\
=&
  \frac{a_{\lambda,j}}{\ln\frac R{s_{\lambda, j}}
}   \bar G(x_{\lambda,i},x_{\lambda,j})+\frac{a_{\lambda,j}}{\ln\frac R{s_{\lambda, j}}
}
\bigl\langle D  \bar G(x_{\lambda,i},x_{\lambda,j}) ,y-x_{\lambda,i}\bigr\rangle +
 O\Bigl( \frac{s_{\lambda,i}^2 }{\ln\lambda
}
\Bigr).
\end{split}
\]
So,  by using \eqref{1-20-10} we obtain the result.

\end{proof}

 For any function $w$, for each $j$,  we denote $\tilde w_j(y)= w(
 s_{\lambda, j} y  + x_{\lambda, j})$.
 In the following, we always assume $L>0$ is a large fixed constant.

\begin{lemma}\label{l1-23-5}
   The set
\[
\Gamma_{\lambda, j}= \bigl\{ y:\;  \mathcal U_{\lambda, {\bf x}_\lambda, {\bf a}_\lambda} ( s_{\lambda, j} y+x_{\lambda, j}  )= \kappa_{\lambda,j}\bigr\}\cap B_{L}(0)
\]
is a closed curve in $\mathbb R^2$, which can be written as
\[
\Gamma_{\lambda, j}(\theta) = (1+ t_{\lambda, j}) ( \cos \theta, \sin\theta),\quad \theta\in [0, 2\pi],
\]
where  $t_{\lambda, j}(\theta)$
is a $C^1$ function satisfying, for some $L>0$,
\begin{equation}\label{15-2-5}
\| t_{\lambda, j}\|_{C^1([0, 2\pi])} \le  Ls_{\lambda, j}.
\end{equation}
  Moreover,
\begin{equation}\label{10-29-5}
\tilde {\mathcal U}_{\lambda, {\bf x}_\lambda, {\bf a}_\lambda,j} ( (1+ t) ( \cos \theta, \sin\theta)  )-\kappa_{\lambda,j}
\begin{cases}
>0, & \text{if}\;\; t<t_{\lambda, j}(\theta);\\
<0,  & \text{if}\;\; t>t_{\lambda, j}(\theta).
\end{cases}
\end{equation}

\end{lemma}

\begin{proof}

It follows from \eqref{1-11-9} that
\begin{equation}\label{12-2-5}
\begin{split}
&\tilde { \mathcal U}_{\lambda, {\bf x}_\lambda, {\bf a}_\lambda,j}(y)-\kappa_{\lambda,i}\\
=&U_{\lambda,x_{\lambda,i}, a_{\lambda,i}}(s_{\lambda, j} y  + x_{\lambda, j})- a_{\lambda,i} -
 \frac {a_{\lambda,i}s_{\lambda,i}}{\ln \frac{R}{s_{\lambda, i}}}
 \bigl\langle D g(x_{\lambda,i},x_{\lambda,i}), y\bigr\rangle
 \\
 &+\sum_{j\ne i}\frac{a_{\lambda,j}s_{\lambda,i}}{\ln\frac R{s_{\lambda, j}}
}
\bigl\langle D  \bar G(x_{\lambda,i},x_{\lambda,j}) ,y\bigr\rangle
 +O\Bigl( \frac{s_{\lambda,i}^2 }{\ln\lambda } \Bigr)
 \\
=&U_{\lambda,x_{\lambda,i}, a_{\lambda,i}}(s_{\lambda, j} y  + x_{\lambda, j})- a_{\lambda,i}+O\Bigl( \frac{s_{\lambda,i} |y| }{\ln\lambda } \Bigr),
\quad y\in B_{L}(0).
\end{split}
\end{equation}

Noting that
\begin{equation}\label{1-16-11}
U_{\lambda,x_{\lambda,i}, a_{\lambda,i}}(s_{\lambda, j} y  + x_{\lambda, j})=
\begin{cases}
a_{\lambda, j}+ \ds\frac{\bar \lambda s^2_{\lambda, j}} {4 } \bigl( 1-|y|^2\bigr),  & |y|\le 1;\vspace{2mm}\\
a_{\lambda, j}\bigl( 1+\ds\frac{\ln |y|}{\ln \frac{s_{\lambda, j}}R  } \bigr), &   |y|\ge 1,
\end{cases}
\end{equation}
by \eqref{1-31-5}, we find  that if
$|y|< 1- L_1 s_{\lambda, j}$ for some large $L_1>0$,  then
\begin{equation}\label{13-2-5}
\begin{split}
\tilde { \mathcal U}_{\lambda, {\bf x}_\lambda, {\bf a}_\lambda,j}(y)-\kappa_{\lambda,i}
=& \frac{\bar \lambda s^2_{\lambda, j}} {4 } \bigl( 1-|y|^2\bigr)+O\Bigl( \frac{s_{\lambda,i} |y| }{\ln\lambda } \Bigr)
\\
> & \frac{\bar \lambda s^2_{\lambda, j}} {4 } \bigl( 1-|1- L_1 s_{\lambda, j}|^2\bigr)+O\Bigl( \frac{s_{\lambda,i} |y| }{\ln\lambda } \Bigr)>0.
\end{split}
\end{equation}

If $ 1+  L_1 s_{\lambda, j}<|y|\le L_2<<L_1$ for some large $L_2$, then
\begin{equation}\label{14-2-5}
\begin{split}
\tilde { \mathcal U}_{\lambda, {\bf x}_\lambda, {\bf a}_\lambda,j}(y)-\kappa_{\lambda,i}
= & \frac{a_{\lambda, j} \ln|y|}{\ln \frac {s_{\lambda, j}}R}
+O\Bigl(\sum_{i=1}^k \frac {L_2 s_{\lambda, i} }{|\ln s_{\lambda, i}|}\Bigr)  \\
< &\frac{a_{\lambda, j} \ln| 1+  L_1 s_{\lambda, j}  |}{\ln \frac {s_{\lambda, j}}R}
+O\Bigl(\sum_{i=1}^k \frac {L_2 s_{\lambda, i} }{|\ln s_{\lambda, i}|}\Bigr)<0.
\end{split}
\end{equation}

Moreover,  it is easy to check that if $L_2<|y|\le \frac\delta{s_{\lambda, j}}$, then
\begin{equation}\label{nnn14-2-5}
\begin{split}
\tilde { \mathcal U}_{\lambda, {\bf x}_\lambda, {\bf a}_\lambda,j}(y)-\kappa_{\lambda,i}
< & U_{\lambda,x_{\lambda,i}, a_{\lambda,i}}(s_{\lambda, j} y  + x_{\lambda, j})- a_{\lambda,i}
 +O\Bigl(\sum_{i=1}^k \frac {1}{|\ln s_{\lambda, i}|}\Bigr) \\
< &\frac{a_{\lambda, j} \ln L_2}{\ln \frac {s_{\lambda, j}}R}
+O\Bigl(\sum_{i=1}^k \frac {1 }{|\ln s_{\lambda, i}|}\Bigr)<0.
\end{split}
\end{equation}
So we have proved that for any $(\cos\theta, \sin\theta)$, there
exists a $t_{\lambda, j}(\theta)$, such that  $|t_{\lambda, j}(\theta)|\le L s_{\lambda, j}$, and
\[
(1+ t_{\lambda, j}(\theta)) ( \cos \theta, \sin\theta)\in \Gamma_{\lambda, j}.
\]

On the other hand, from \eqref{1-11-9}, \eqref{1-6-11}  and \eqref{1-16-11},  we have
\begin{equation}\label{30-23-5}
\begin{split}
& \frac{\partial \tilde { \mathcal U}_{\lambda, {\bf x}_\lambda, {\bf a}_\lambda,j}( (1+ t) ( \cos \theta, \sin\theta)  )}{\partial t}
\Bigr|_{t=0}\\
=&
\frac{\partial U_{\lambda,0,a_{\lambda, j}}( (1+ t) ( \cos \theta, \sin\theta))}
{\partial t}\Bigr|_{t=0}\\
&
-
 \frac {a_{\lambda,j}}{\ln \frac{R}{s_{\lambda, j}}}
 \bigl\langle D g(x_{\lambda, j},x_{\lambda, j}), ( \cos \theta, \sin\theta))\bigr\rangle s_{\lambda,j}
 \\
 &+\sum_{i\ne j}\frac{a_{\lambda,i}}{\ln\frac R{s_{\lambda, i}}
}
\bigl\langle D  \bar G(x_{\lambda, i},x_{\lambda, j}) ,( \cos \theta, \sin\theta))\bigr\rangle s_{\lambda,j}+O\Bigl(\sum_{i=1}^k\frac {s_{\lambda,j}^2}{|\ln s_{\lambda, i}|}\Bigr)\\
=& -\frac{\lambda s_{\lambda, j}^2}{2}+O\Bigl( \frac{|\nabla \mathcal W ({\bf x}_\lambda)|}{\ln \lambda} s_{\lambda,j}+\sum_{i=1}^k\frac {s^2_{\lambda,j}}{|\ln s_{\lambda, i}|^2}\Bigr)<  0.
\end{split}
\end{equation}
    So $t_{\lambda,j}$ is unique.

    Differentiating   $\tilde {\mathcal U}_{\lambda, {\bf x}_\lambda, {\bf a}_\lambda}( (1+ t_{\lambda, j}(\theta)) ( \cos \theta, \sin\theta)=\kappa_{\lambda,j} $ with respect to $\theta$,  noting that
\[
    \begin{split}
    &\frac{\partial \widetilde {PU}_{\lambda,x_{\lambda,  i},a_{\lambda, i}, j}(1+ t_{\lambda, j}(\theta)) ( \cos \theta, \sin\theta)}{\partial \theta}\\
    =&O\Bigl( \frac{s_{\lambda,j}}{\ln \lambda}
\Bigr)\bigl| \frac{\partial t_{\lambda,j}}{\partial \theta}\bigr|+O\Bigl(\sum_{i=1}^k \frac {s_{\lambda,i}}{|\ln s_{\lambda,i}|}\Bigr), \quad  \forall\; i\ne j,
    \end{split}
    \]
we find
\begin{equation}\label{1-29-5}
\begin{split}
&
 \frac{\partial   \widetilde {PU}_{\lambda,x_{\lambda, j},a_{\lambda, j}, j}
  ( (1+ t_{\lambda,j}(\theta)) ( \cos \theta, \sin\theta))}
{\partial \theta}
+O\Bigl( \frac{s_{\lambda,j}}{\ln \lambda}
\Bigr)
\bigl| \frac{\partial t_{\lambda,j}}{\partial \theta}\bigr|\\
= & O\Bigl(\sum_{i=1}^k \frac {s_{\lambda,i}}{|\ln s_{\lambda,i}|}\Bigr).
\end{split}
\end{equation}
Similar to    \eqref{30-23-5},  we can estimate
\[
\frac{\partial   \widetilde {PU}_{\lambda,x_{\lambda, j},a_{\lambda, j}, j}
  ( (1+ t_{\lambda,j}(\theta)) ( \cos \theta, \sin\theta))}
{\partial \theta}=
\Bigl[-\frac{\lambda s_{\lambda, j}^2}{2}+O\Bigl( \frac{s_{\lambda,j}}{\ln \lambda} \Bigr)\Bigr] \frac{\partial t_{\lambda,j}}{\partial \theta},
\]
which, together with \eqref{1-29-5},
gives
\[
 \frac{\partial t_{\lambda, j}}{\partial \theta}= O(s_{\lambda, j}).
\]
So $\Gamma_{\lambda, j}$ is a smooth closed curve in $\mathbb R^2$
and
\eqref{15-2-5} satisfies. It is also easy to check that \eqref{10-29-5}
holds.

\end{proof}

\begin{lemma}\label{al1}

 The set
\[
\tilde \Gamma_{\lambda, j}= \bigl\{ y:\;  u_\lambda( s_{\lambda, j }y+ x_{\lambda, j})  = \kappa_{\lambda,j}\bigr\}\cap B_{L}(0)
\]
is a continuous closed curve in $\mathbb R^2$, and
\begin{equation}\label{12-29-5}
\begin{split}
&\tilde \Gamma_{\lambda, j}(\theta) = (1+\tilde t_{\lambda, j}(\theta) ) ( \cos \theta, \sin\theta)\\
=&\Bigl(1+\frac1{a_{\lambda, j}}\tilde\omega_{\lambda, j}( (1+\tilde t_{\lambda, j}(\theta) ) \cos \theta, (1+\tilde t_{\lambda, j}(\theta) )
\sin\theta)\ln\frac R{s_{\lambda, j}} \Bigr) \bigl( \cos \theta, \sin\theta\bigr)\\
&- s_{\lambda, j}  \bigl\langle
\nabla g(x_{\lambda, j},x_{\lambda, j}
)-\sum_{i\ne j}\frac{ a_{\lambda, i}}{a_{\lambda, j}}\nabla \bar G(x_{\lambda, i},x_{\lambda, j}) , ( \cos \theta, \sin\theta) \bigr\rangle \bigl( \cos \theta, \sin\theta\bigr)\\
&+ O\Bigl( \frac1\lambda +|\tilde \omega_{\lambda, j} ( (1+\tilde t_{\lambda, j}(\theta) ) \cos \theta, (1+\tilde t_{\lambda, j}(\theta) ) \sin\theta)|^2|\ln s_{\lambda, j}|^2\Bigr),\quad \theta\in [0, 2\pi],
\end{split}
\end{equation}
for some function $\tilde t_{\lambda, j}(\theta)$.
 Moreover,
\begin{equation}\label{1-18-8}
|\tilde \omega_{\lambda, j} ( (1+\tilde t_{\lambda, j}(\theta) ) \cos \theta, (1+\tilde t_{\lambda, j}(\theta) ) \sin\theta)||\ln s_{\lambda, j}|\to 0,
\end{equation}
\begin{equation}\label{11-29-5}
u_\lambda( (1+ t) ( \cos \theta, \sin\theta)  )-\kappa_{\lambda,j}
\begin{cases}
>0, & \text{if}\; t<\tilde t_{\lambda, j}(\theta);\\
<0,  & \text{if}\; t>\tilde t_{\lambda, j}(\theta).
\end{cases}
\end{equation}

\end{lemma}

\begin{proof}

  If follows from  Lemma~\ref{p1-3-8} that
\begin{equation}\label{31-28-5}
 \frac{\partial \tilde  u_{\lambda} ( (1+t)(\cos\theta, \sin\theta) )}{\partial t}
\Bigr|_{t=0}\\
=\bar \lambda r_{\lambda, j}^2\Bigl( \frac1{r_{\lambda, j}}(w'(1)+o(1))\Bigr) <  0.
\end{equation}
As a result,  $\tilde t_{\lambda, j}$  is unique.   Therefore,  $\Gamma_{\lambda, j}$
is a continuous closed curve in $\mathbb R^2$. Moreover, Lemma~\ref{p1-3-8}
also implies that for any $y_\lambda \in \tilde \Gamma_{\lambda, j}$, $|y_\lambda|\to 1$ as $\lambda\to +\infty$.

It is easy to check from \eqref{31-28-5} that \eqref{11-29-5} holds.

 For any point
\[
y_{\lambda, j}(\theta)= (1+ \tilde t_{\lambda, j}(\theta) )(\cos\theta, \sin\theta) \in \tilde \Gamma_{\lambda, j},
\]
it follows from \eqref{1-11-9} that if $|y_{\lambda,j}(\theta)|\ge 1$, then
\begin{equation}\label{32-28-5}
\begin{split}
0=&\tilde {\mathcal U}_{\lambda, {\bf x}_\lambda, {\bf a}_\lambda,j}( y_{\lambda,j}(\theta)  )-\kappa_{\lambda,j}+\tilde\omega_{\lambda, j}(y_{\lambda,j}(\theta))\\
=&
 U_{\lambda, x_{\lambda,j}, a_{\lambda,j}, j}( s_{\lambda, j} y + x_j )- a_{\lambda,j} -
 \frac {a_{\lambda,j}s_{\lambda, j}}{\ln \frac{R}{s_{\lambda, j}}}
 \bigl\langle D g(x_j,x_j), y_{\lambda,j}\bigr\rangle
 \\
 &+\sum_{i\ne j}\frac{a_{\lambda,i}s_{\lambda, j}}{\ln\frac R{s_{\lambda, i}}
}
\bigl\langle D  \bar G(x_{\lambda, i}, x_{\lambda, j}) , y_{\lambda,j}\bigr\rangle
 +O\Bigl( \frac{1 }{\lambda \ln\lambda } \Bigr)+\tilde\omega_{\lambda, j}(y_{\lambda,j}(\theta))\\
=& a_{\lambda, j} \frac{\ln|y_{\lambda, j}|}{\ln\frac{s_{\lambda, j}}R}+\tilde\omega_{\lambda, j}(y_{\lambda, j})
 -
 \frac {a_{\lambda,j}s_{\lambda, j}}{\ln \frac{R}{s_{\lambda, j}}}
 \bigl\langle D g(x_{\lambda, j},x_{\lambda, j}), y_{\lambda,j}\bigr\rangle
 \\
 &+\sum_{i\ne j}\frac{a_{\lambda,i}s_{\lambda, j}}{\ln\frac R{s_{\lambda, i}}
}
\bigl\langle D  \bar G(x_{\lambda, i},x_{\lambda, j}) , y_{\lambda,j}\bigr\rangle
 +O\Bigl( \frac{1 }{\lambda \ln  \lambda } \Bigr).
\end{split}
\end{equation}
So
\begin{equation}\label{33-28-5}
\begin{split}
|y_{\lambda, j}| =& e^{ \frac{ \bigl(\tilde\omega_{\lambda, j}(y_{\lambda, j})
 -
 \frac {a_{\lambda,j}s_{\lambda, j}}{\ln \frac{R}{s_{\lambda, j}}}
 \bigl\langle D g(x_{\lambda, j},x_{\lambda, j}), y_{\lambda,j}\bigr\rangle
 +\sum_{i\ne j}\frac{a_{\lambda,i}s_{\lambda, j}}{\ln\frac R{s_{\lambda, i}}
}
\bigl\langle D  \bar G(x_{\lambda, i},x_{\lambda, j}) , y_{\lambda,j}\bigr\rangle
 +O( \frac{1 }{\lambda \ln  \lambda } )\bigr)\ln
\frac R{s_{\lambda, j}} }{a_{\lambda, j}}}
\\
=&1+\tilde\omega_{\lambda, j}(y_{\lambda, j})\frac {\ln \frac{R}{s_{\lambda, j}}}{a_{\lambda,j}}
 -s_{\lambda, j}
 \langle D g(x_{\lambda, j},x_{\lambda, j}), y_{\lambda,j}\rangle
 +\sum_{i\ne j}\frac{a_{\lambda,i}s_{\lambda, j}}{a_{\lambda,j}
}
\langle D  \bar G(x_{\lambda, i},x_{\lambda, j}) , y_{\lambda,j}\rangle\\&
 +  O\Bigl( \frac1\lambda+(\tilde\omega_{\lambda, j}(y_{\lambda, j}) \ln s_{\lambda, j})^2\Bigr).
\end{split}
\end{equation}

If $|y_{\lambda,j}(\theta)|< 1$, then
\begin{equation}\label{34-28-5}
\begin{split}
0=&\tilde {\mathcal U}_{\lambda, {\bf x}_\lambda, {\bf a}_\lambda,j}( y_{\lambda,j}(\theta)  )-\kappa_{\lambda,j}+\tilde\omega_{\lambda, j}(y_{\lambda,j}(\theta))
\\
=& \frac{\bar \lambda s_{\lambda,j}^2}{4} \bigl( 1-|y_{\lambda, j}|^2\bigr)+\tilde\omega_{\lambda, j}(y_{\lambda, j})
 -
 \frac {a_{\lambda,j}s_{\lambda, j}}{\ln \frac{R}{s_{\lambda, j}}}
 \bigl\langle D g(x_{\lambda, j},x_{\lambda, j}), y_{\lambda,j}\bigr\rangle
 \\
 &+\sum_{i\ne j}\frac{a_{\lambda,i}s_{\lambda, j}}{\ln\frac R{s_{\lambda, i}}
}
\bigl\langle D  \bar G(x_{\lambda, i},x_{\lambda, j}) , y_{\lambda,j}\bigr\rangle
 +O\Bigl( \frac{1}{\lambda \ln \lambda } \Bigr).
\end{split}
\end{equation}
So,
\begin{equation}\label{35-28-5}
\begin{split}
|y_{\lambda, j}| =& 1+ \frac {2}{\bar \lambda s_{\lambda, j}^2} \Bigl(
\tilde\omega_{\lambda, j}(y_{\lambda, j})
 -
 \frac {a_{\lambda,j}s_{\lambda, j}}{\ln \frac{R}{s_{\lambda, j}}}
 \bigl\langle D g(x_{\lambda, j},x_{\lambda, j}), y_{\lambda,j}\bigr\rangle
 +\sum_{i\ne j}\frac{a_{\lambda,i}s_{\lambda, j}}{\ln\frac R{s_{\lambda, i}}
}
\bigl\langle D  \bar G(x_{\lambda, i},x_{\lambda, j}) , y_{\lambda,j}\bigr\rangle
\Bigr)\\
& +O\Bigl( \frac1\lambda+(\tilde\omega_{\lambda, j}(y_{\lambda, j}) \ln s_{\lambda, j})^2\Bigr)\\
=&1+\tilde\omega(y_{\lambda, j})\frac {\ln \frac{R}{s_{\lambda, j}}}{a_{\lambda,j}}
 -s_{\lambda, j}
 \langle D g(x_{\lambda,j},x_{\lambda, j}), y_{\lambda,j}\rangle
 +\sum_{i\ne j}\frac{a_{\lambda,i}s_{\lambda, j}}{a_{\lambda,j}
}
\langle D  \bar G(x_{\lambda, i},x_{\lambda, j}) , y_{\lambda,j}\rangle
\\
& + O\Bigl( \frac1\lambda+(\tilde\omega_{\lambda, j}(y_{\lambda, j}) \ln s_{\lambda, j})^2\Bigr).
\end{split}
\end{equation}
So we find  \eqref{12-29-5}  follows from \eqref{33-28-5}  and \eqref{35-28-5}.
On the other hand,
 it is easy to see that  \eqref{32-28-5}, \eqref{34-28-5}  and $|y_{\lambda, j}|\to 1$  implies \eqref{1-18-8}.

\end{proof}

\end{document}